\documentclass[10pt]{article}
\usepackage{graphicx}
\usepackage{epstopdf}
\usepackage{amsmath,amssymb,amsthm,amsfonts}
\usepackage{color}
\usepackage{hyperref}
\usepackage{cite}
\usepackage{verbatim}
\usepackage{appendix}

\allowdisplaybreaks[4]

\newtheorem{thm}{Theorem}[section]

\newtheorem{lem}[thm]{Lemma}
\newtheorem{prop}[thm]{Proposition}

\theoremstyle{definition}

\theoremstyle{remark}
\newtheorem{rem}{Remark}[section]

\numberwithin{equation}{section}

\DeclareMathSymbol{\C}{\mathalpha}{AMSb}{"43}

\newcommand{\eps}{\varepsilon}
\newcommand{\Om}{\Omega}

\newcommand{\gam}{\gamma}
\newcommand{\B}{{\beta}}
\newcommand{\R}{{\mathbb{R}}}
\newcommand{\h}{{\mathcal{H}}}
\newcommand{\m}{{\mathcal{M}}}

\newcommand{\inte}{\int_{\mathbb{R}^2}}

\newcommand{\f}{a_1,a_2,\beta}

\def\ds{\displaystyle}

\def\R{{\mathbb R}}
\def\C{{\mathbb C}}
\def\ovl{\overline}
\def\l{\left}
\def\r{\right}
\def\p{\partial}

\newcommand{\bsub}{\begin{subequations}}
\newcommand{\esub}{\end{subequations}$\!$}
\title {Axial Symmetry and Refined Spike Profiles of Ground States for Rotating Two-Component Bose Gases}
\author{ Yongshuai Gao\thanks{Email: ysgao@mails.ccnu.edu.cn. Y. Gao is partially supported by the Graduate Education Innovation Funds $\#$2022CXZZ088 at Central China Normal University in P. R. China.},
and Yong Luo\thanks{Email: yluo@ccnu.edu.cn. Y. Luo is partially supported by NSFC under Grant 12201231.}
\\
\small \it	School of Mathematics and Statistics,\\
\small \it  Hubei key Laboratory of Mathematical Sciences,\\
\small \it Central China Normal University, Wuhan 430079, P. R. China\\}

\begin{document}
\baselineskip= 15pt
\maketitle
\begin{abstract}
This paper considers ground states of two-component Bose gases confined in an anharmonic trap 
rotating at the velocity $0<\Om<\infty$, where the intraspecies interaction $(-a_1,-a_2)$ and the interspecies interaction $-\B$ are both attractive, i.e, $a_1, a_2$ and $\B $ are all positive. We prove the axially symmetry and
the refined spike profiles of ground states as $\B\nearrow\B^*:=a^*+\sqrt{(a^*-a_1)(a^*-a_2)}$, where $0<a_1,\, a_2<a^*:=\|Q\|^2_2$ are fixed and $Q>0$ is the unique positive solution of $-\Delta u+u-u^3=0$ in $\R^2$.



\end{abstract}

\section{Introduction}
In this paper, we study ground states of the following two coupled rotational Gross-Pitaevskii equations:
\begin{equation}\label{1.1:GPs}
\left\{\begin{array}{lll}
-\Delta u_1+V(x)u_1+i\, \Om \big(x^{\perp}\cdot \nabla u_1\big)
=\mu u_1+a_1|u_1|^2u_1+\B |u_2|^2u_1\  \ \hbox{in}\  \ \R^2,\\[3mm]
-\Delta u_2+V(x)u_2+i\, \Om \big(x^{\perp}\cdot \nabla u_2\big)
=\mu u_2+a_2|u_2|^2u_2+\B |u_1|^2u_2\  \ \hbox{in}\  \ \R^2,\\[3mm]
\ds\inte (|u_{1}|^2+|u_{2}|^2)dx=1,
\end{array}\right.
\end{equation}
where $x^\perp=(-x_2,x_1)$ holds for $x=(x_1,x_2)\in\R^2$, and
$\mu=\mu (a_1, a_2, \B)\in\R$ is the chemical potential.
The system \eqref{1.1:GPs} is used to model the rotating two-component Bose gases \cite{GRPG,KTU1,LS,LS1}, where $V(x)\geq0$ is a trapping potential rotating at the velocity $\Om>0$, $a_j>0$ (resp. $<0$) represents
the attractive (resp. repulsive) intraspecies interaction of cold atoms inside the $j^{th}$ component ($j=1,2$), and while $\B>0$ (resp. $<0$) describes the attractive (resp. repulsive) interspecies interaction of cold atoms between two components.

When there is no rotation for the system, i.e. $\Om=0$, \eqref{1.1:GPs} can be reduced to the real-valued problem.
In this case, there are many rich results, including the existence and nonexistence, the limiting behavior, symmetry breaking and the local uniqueness of spike solutions for \eqref{1.1:GPs}, see \cite{BC,GLWZ1,GX,LW,R,ZZ,WY} and the references therein.
In particular, the authors recently proved in \cite{ZZ} the single speak solutions of \eqref{1.1:GPs} are cylindrically symmetric and unique for the ring-shaped potential $V(x)=(|x|-1)^2$.

When there is rotation for the system, i.e. $\Om>0$, ground states of \eqref{1.1:GPs} are no longer real-valued. For the following harmonic trapping potential:
\begin{equation*}
V(x)=\lambda_1^2x_1^2+\lambda_2^2x_2^2, \ \ \hbox{where} \ \ \lambda_1,\lambda_2>0,
\end{equation*}
there exist some interesting works \cite{A,AMW,ANS,AS,GGLL} on the existence, the limiting behavior and the nonexistence of vortices for ground states.
Specially, 
by developing the argument of refined expansions, we recently proved in \cite{GGLL} that the nonexistence of vortices for ground states in the two-component focusing (i.e., $a_1,a_2,\B>0$) Bose gases. We remark that the above works mainly focused on the case $0<\Om<\Om^*$, where $\Om^*:=2\min\{\lambda_1,\lambda_2\}$ is the critical rotational speed.

In order to investigate the rotating Bose-Einstein condensates at arbitrarily large rotational speed $\Om>0$,
physicists usually replace the harmonic potential by the following anharmonic potential:
\begin{equation}\label{1.3:quar}
V(x)=|x|^4+\Lambda_1x_1^2+\Lambda_2x_2^2, \ \ \hbox{where} \ \ \Lambda_1,\Lambda_2\in \R.
\end{equation}
As for the rotating Bose-Einstein condensates confined in the potential \eqref{1.3:quar}, we refer to \cite{BSSD,F,F2,SBCD,KMKR,LZ,TKT} and the references therein.

Motivated by the above works mentioned \cite{GGLL,GLLP,GLP,GLW,GLY}, in the present paper we study ground states of
\eqref{1.1:GPs}, where the intraspecies interaction and interspecies interaction are both attractive, i.e., $a_1, a_2, \B>0$, and $V(x)\ge0$ is an anharmonic trapping potential satisfying
\begin{equation}\label{1.4:Vx}
V(x)=|x|^4+\frac{\Om^2-8}{4}|x|^2+1, \ \ \hbox{where} \ \ 0<\Om<\infty
\ \ \hbox{is fixed}.
\end{equation}
Under the anharmonic trapping potential \eqref{1.4:Vx}, the system \eqref{1.1:GPs} can be rewritten as the following form
\begin{equation}\label{1.5:GPs}
\left\{\begin{array}{lll}
\ds-\Big(\nabla-i\frac{\Om}{2}x^\perp\Big)^2u_1+\big(|x|^2-1\big)^2u_1
=\mu u_1+a_1|u_1|^2u_1+\B |u_2|^2u_1\  \ \hbox{in}\  \ \R^2,\\[3mm]
\ds-\Big(\nabla-i\frac{\Om}{2}x^\perp\Big)^2u_2+\big(|x|^2-1\big)^2u_2
=\mu u_2+a_2|u_2|^2u_2+\B |u_1|^2u_2\  \ \hbox{in}\  \ \R^2,\\[3mm]
\ds\inte (|u_{1}|^2+|u_{2}|^2)dx=1,
\end{array}\right.
\end{equation}
where $-\big(\nabla-i\frac{\Om}{2}x^\perp\big)^2
=-\Delta+i\, \Om x^{\perp}\cdot \nabla+\frac{\Om^2}{4}|x|^2$.
When $\Om=0$ and the ring-shaped potential $\big(|x|^2-1\big)^2$ is replaced by some potentials $V(x)$ with some non-degenerate critical points, the authors in \cite{GX,GLWZ1} investigated the existence and uniqueness of the spike solutions for \eqref{1.5:GPs}.
One can note that the ring-shaped potential $\big(|x|^2-1\big)^2$ in \eqref{1.5:GPs} has infinitely critical points on the sphere $|x|=1$ and its directional derivative along the tangent direction is always zero at each critical point lying on the sphere $|x|=1$.
Thus, the critical points of the ring-shaped potential $\big(|x|^2-1\big)^2$ do not satisfy the non-degenerate conditions of \cite{GX,GLWZ1}. Therefore, the arguments of \cite{GX,GLWZ1} is not applicable to the analysis of \eqref{1.5:GPs}.


Similar to \cite [Prosition A.1] {GLWZ1},  ground states of
\eqref{1.1:GPs} under the anharmonic trapping potential \eqref{1.4:Vx} can be described equivalently by minimizers of the following complex-valued constraint variational problem:
\begin{equation}\label{1.6:cvp}
e(\Om,\f):=\inf_{(u_1,u_2)\in\m}F_{\Om,\f}(u_1,u_2),
\  \ \Om >0,\ \, a_1>0,\  \ a_2>0,\  \ \B>0,
\end{equation}
where the  energy functional $F_{\Om,\f}(u_1,u_2)$ is given by
\begin{equation}\label{1.7:GPf}
\begin{split}
 F_{\Om,\f}(u_1,u_2):=
 &\sum_{j=1}^2\inte\Big[\Big|\Big(\nabla-i\frac{\Om}{2}x^\perp\Big) u_j\Big|^2
 +\big(|x|^2-1\big)^2|u_j|^2
 -\frac{a_j}{2}|u_j|^4\Big]dx\\
  &
 -\inte\B|u_1|^2|u_2|^2dx,
\end{split}
\end{equation}
and the space $\m$ is defined as
\begin{equation}\label{1.8:m}
\begin{split}
\m:=\Big\{(u_1,u_2)\in \h \times \h:
\ \ \inte (|u_{1}|^2+|u_{2}|^2)dx=1\Big\}.
\end{split}
\end{equation}
Here $(iu_j,\nabla u_j)=i(u_j\nabla \bar u_j-\bar u_j\nabla u_j)/2$, and
\begin{equation*}
\h:=\Big\{u\in  H^1(\R^2,\C):\ \inte \big(|x|^2-1\big)^2|u(x)|^2dx<\infty \Big\}.
\end{equation*}
In order to study ground states of \eqref{1.1:GPs}, we shall therefore focus on minimizers of \eqref{1.6:cvp} in this paper.

Throughout the whole paper, we always use $Q=Q(|x|)>0$ (cf. \cite{GNN,K}) to denote the unique positive solution of the following equation
\begin{equation}\label{1.10:Q}
-\Delta u+ u-u^3=0\  \  \mbox{in} \ \ \R^2,  \ \ u\in H^1(\R^2,\R),
\end{equation}
and denote $\B^*>0$ by
\begin{equation}\label{1.11:B*}
\B^*:=\B^*(a_1,a_2)=a^*+\sqrt{(a^*-a_1)(a^*-a_2)},\ \  0<a_1,a_2<a^*:=\|Q\|_2^2.
\end{equation}

When the trapping potential $V(x)\ge0$ satisfies \eqref{1.4:Vx},
it follows immediately from \cite [Theorems 2.1 and 2.2] {GGLL} the following existence and nonexistence of minimizers for $e(\Om,\f)$:

\noindent{\bf Theorem A.1} {\em(Theorems 2.1 and 2.2 in \cite{GGLL})
Suppose $V(x)$ satisfies \eqref{1.4:Vx}. Then for any given $0<\Om<\infty$, $e(\Om,\f)$ admits minimizers if and only if $0<a_1,a_2<a^*$ and $0<\B<\B^*$.}

Let $(u_{1\B},u_{2\B})$ be a minimizer of $e(\Om,\f)$, we then obtain from the variational theory that $(u_{1\B},u_{2\B})$ solves the following Euler-Lagrange system
\begin{equation}\label{1.12:ELs}
\left\{\begin{array}{lll}
\ds-\Big(\nabla-i\frac{\Om}{2}x^\perp\Big)^2u_{1\B}+\big(|x|^2-1\big)^2u_{1\B}
=\mu_\B u_{1\B}+a_1|u_{1\B}|^2u_{1\B}+\B |u_{2\B}|^2u_{1\B}\  \ \hbox{in}\  \ \R^2,\\[3mm]
\ds-\Big(\nabla-i\frac{\Om}{2}x^\perp\Big)^2u_{2\B}+\big(|x|^2-1\big)^2u_{2\B}
=\mu_\B u_{2\B}+a_2|u_{2\B}|^2u_{2\B}+\B |u_{1\B}|^2u_{2\B}\  \ \hbox{in}\  \ \R^2,\\[3mm]
\ds\inte (|u_{1\B}|^2+|u_{2\B}|^2)dx=1,
\end{array}\right.
\end{equation}
where $\mu_\B\in\R$ is a suitable Lagrange multiplier and satisfies
\begin{equation}\label{1.13:mu}
\mu_\B=e(\Om,\f)-\inte \Big(\frac{a_1}{2}|u_{1\B}|^4+\frac{a_2}{2}|u_{2\B}|^4+\B|u_{1\B}|^2|u_{2\B}|^2\Big)dx.
\end{equation}
Similar to the arguments of (3.12) and (3.13) in \cite {GGLL}, one can deduce that
\begin{equation}\label{1.14:mu-infty}
\mu_\B\to-\infty \ \ \hbox{and}\ \  \|u_{j\beta}\|_{L^\infty(\R^2)}\to\infty\ \ (j=1,2)\ \  \ \ \hbox{as} \ \ \B\nearrow\B^*.
\end{equation}
The main purpose of this paper is to investigate {\em the axial symmetry and the refined spike profiles of minimizers for \eqref{1.6:cvp} as $\B\nearrow\B^*$.}

\subsection{Main results}
In this subsection, we shall introduce the main results of the present paper.
For convenience, we denote $\gam_j>0$ and $\alpha_\B>0$ by
\begin{equation}\label{1.15:gam}
\gam_j:=1-\frac{\sqrt{a^*-a_j}}{\sqrt{a^*-a_1}+\sqrt{a^*-a_2}}\in(0,1),
\ \  0<a_1,a_2<a^*, \ \ j=1,2,
\end{equation}
and
\begin{equation}\label{1.16:alp}
\alpha_\B:=\Big(\frac{8\gamma_1\gamma_2(\B^*-\B)}{(\Om^2+8)\inte|x|^2Q^2}\Big)^\frac{1}{4}>0,
\ \ 0<\B<\B^*.
\end{equation}

Using the above notations, our first result is to prove the following limiting behavior of minimizers for $e(\Om,\f)$ as $\B\nearrow\B^*$.

\begin{thm}\label{thm1}
Suppose $V(x)$ satisfies \eqref{1.4:Vx}, and assume $0<\Om<\infty$. Let $(u_{1\B},u_{2\B})$ be a complex-valued minimizer of $e(\Om,\f)$, where $0< a_1,a_2<a^*$ and $0<\B<\B^*$. Then for $j=1,2$, we have
\begin{equation}\label{1.17:lim}
\widetilde{v}_{j\B}(x):=\sqrt{a^*}\alpha_\B u_{j\B}(\alpha_\B x+x_\B)
e^{-i\big(\frac{\alpha_\B\Om}{2}x\cdot x_\B^\perp-\widetilde{\theta}_{j\B}\big)}
\to\sqrt{\gamma_j}Q(x)\ \ \hbox{as}\ \ \B\nearrow\B^*
\end{equation}
strongly in $H^1(\R^2,\C)\cap L^\infty(\R^2,\C)$,
where $\alpha_\B>0$ is given by \eqref{1.16:alp}, $\gam_j\in(0,1)$ is as in \eqref{1.15:gam}, $\widetilde{\theta}_{j\B}\in [0,2\pi)$ is a properly chosen constant, and $x_\B\in\R^2$ is the unique maximal point of $|u_{1\B}|^2+|u_{2\B}|^2$ and satisfies
\begin{equation}\label{1.18:xB}
x_\B\to P_0\ \ \hbox{as}\ \ \B\nearrow\B^*
\end{equation}
for some $P_0\in\R^2$ satisfying $|P_0|=1$.
\end{thm}

Theorem \ref{thm1} tells us that the minimizers concentrate at a minimal point $P_0$ of the potential $(|x|^2-1)^2$. Due to the degeneracy of the potential $(|x|^2-1)^2$, it is difficult to obtain the blow up rate by applying the energy estimate methods used in \cite{GLWZ1,GZZ1}.
To overcome this difficulty, we shall prove \eqref{1.17:lim} by applying Pohozaev identities. 
We remark that our method is more simpler and direct than the energy estimate methods used in \cite{GLWZ1,GZZ1}.

Based on Theorem \ref{thm1}, the second result of this paper is the axial symmetry of minimizers for $e(\Om,\f)$, which is stated as the following theorem.

\begin{thm}\label{thm2}
Suppose $V(x)$ satisfies \eqref{1.4:Vx}, and assume $0<\Om<\infty$ and
$0< a_1,a_2<a^*$ are fixed. Let $(u_{1\B},u_{2\B})$ be a complex-valued minimizer of $e(\Om,\f)$ and $x_\B\in\R^2$ be the unique maximal point of $|u_{1\B}|^2+|u_{2\B}|^2$. Then up to the constant phase, the minimizer $(u_{1\B},u_{2\B})$ must be axially symmetric with respect to the $Ox_\B$-line when $\B^*-\B>0$ is small enough.
\end{thm}

To the best of our knowledge, there are few works prove the symmetry of complex-valued solutions about the equation involved magnetic Laplacian, see \cite{BNV,BDELL}.  However, their methods is not applicable to deal with ring-shaped potentials since it is not monotonic.  Moreover, it is hard to use the symmetry argument about \eqref{1.1:GPs} with $\Omega=0$ in \cite{ZZ} to deal with our case $\Omega>0$. As far as we know, this seems the first work on the axial symmetry of solutions for the system of GP equations.

\begin{rem}
As a byproduct, we obtain the local uniqueness of minimizers in the sense that up to the constant phase and the rotational transformation, there is a unique minimizer $(u_{1\B},u_{2\B})$ of $e(\Om,\f)$ when $\B^*-\B>0$ is small enough, see details in Proposition \ref{prop3.1}.
 Without loss of generality, we may assume that the unique maximal point $x_\B$ of $|u_{1\B}|^2+|u_{2\B}|^2$ lies on the $x_2$-axis, i.e., $x_\B=(0,p_\B)$, then the axially symmetric in Theorem \ref{thm2} means that
\begin{equation*}
\big(u_{1\B}(-x_1,x_2),u_{2\B}(-x_1,x_2)\big)\equiv
\big(u_{1\B}(x)e^{i\phi_{1\B}},
     u_{2\B}(x)e^{i\phi_{2\B}}\big)
\end{equation*}
for some constant phase $(\phi_{1\B},\phi_{2\B})\in[0,2\pi)\times[0,2\pi)$ when $\B^*-\B>0$ is small enough.
\end{rem}


{\em Main idea of proving Theorem \ref{thm2}}. The proof of Theorem \ref{thm2} is divided into two steps.

In the first step, we shall prove in Proposition \ref{prop3.1} that up to the constant phase and rotational transformation, the minimizers of $e(\Om,\f)$ must be unique
when $\B^*-\B>0$ is small enough.
Motivated by \cite{GGLL,GLLP,GLP,GLW}, we shall achieve it by
constructing various Pohozaev identities and employing the non-degeneracy of \eqref{1.10:Q}.
Compared with \cite{G,LPWY,ZZ}, since $(|x|^2-1)^2$ is degenerate,
one cannot obtain the useful information from the Pohozaev identities corresponding to the tangent direction of the sphere $|x|=1$. We shall overcome this difficulty by making full use of the rotational invariance of the system \eqref{1.1:GPs} and the new transformation \eqref{3.7}.

Let $(u_{1\B},u_{2\B})$ be a minimizer of $e(\Om,\f)$. Without loss of generality, we may assume the unique maximal point $x_\B$ of $|u_{1\B}|^2+|u_{2\B}|^2$ lies on the $x_2$-axis, i.e., $x_\B=(0,p_\B)$. In the second step, we shall prove that $(u_{1\B}(-x_1,x_2),u_{2\B}(-x_1,x_2))$ is also a minimizer of $e(\Om,\f)$.
Following Proposition \ref{prop3.1}, we finally obtain that
\begin{equation*}
\big(u_{1\B}(-x_1,x_2),u_{2\B}(-x_1,x_2)\big)\equiv
\big(u_{1\B}(x)e^{i\phi_{1\B}},
     u_{2\B}(x)e^{i\phi_{2\B}}\big)
\end{equation*}
for some constant phase $(\phi_{1\B},\phi_{2\B})\in[0,2\pi)\times[0,2\pi)$ when $\B^*-\B>0$ is small enough.
This implies that Theorem \ref{thm2} holds true.

The third result of this paper is the refined spike profiles of minimizers for $e(\Om,\f)$ in terms of $\alpha_\B$ as $\B\nearrow\B^*$, where $\alpha_\B>0$ is as in \eqref{1.16:alp}.
For convenience, we denote $\psi_{1}(x)\in C^2(\R^2,\C)\cap L^\infty(\R^2,\C)$ to be the unique solution of the following equation
\begin{equation}\label{1.19:psi1}
\left\{\begin{array}{lll}
\begin{split}
-\Delta\psi_{1}(x)+\psi_{1}(x)-Q^2\psi_{1}(x)-2Re(Q\psi_{1}(x))Q
             &=-\big(\frac{\Om^2}{4}|x|^2+4x_{2}^2\big)Q
             \,\ \mbox{in}\,\ \R^2,\\
\nabla Re(\psi_{1}(0))=0,\ Re\inte \psi_1(iQ)&=0,
\end{split}
\end{array}\right.
\end{equation}
where $Q>0$ is the unique positive solution of \eqref{1.10:Q}. We also denote
$C_1(x)$ and $C_2(x)$ by
\begin{equation}\label{1.20:Cx}
\begin{split}
C_{1}(x):=&\frac{B_2}{4B_1}\big(Q+x\cdot\nabla Q\big),\\[3mm]
C_{2}(x):=&\frac{8B_1^2B_4-7B_2^2}{32B_1^2}\big(Q+x\cdot\nabla Q\big)
             +\frac{B_2^2}{32B_1^2}x^T(\nabla^2 Q)x
          +\frac{B_2^2}{16B_1^2}\big(x\cdot\nabla Q\big)+\psi_{1}(x),
\end{split}
\end{equation}
where $Q>0$ is the unique positive solution of \eqref{1.10:Q}, $\psi_{1}(x)\in C^2(\R^2,\C)\cap L^\infty(\R^2,\C)$ is the unique solution of \eqref{1.19:psi1},
and $B_1,B_2,B_3,B_4$ are defined by
\begin{equation}\label{1.21:B1234}
\begin{split}
B_{1}:=\frac{\Om^2+8}{4}\inte|x|^2Q^2, \quad
B_2:=2\inte\big(|x|^4+4A|x|^2\big)Q^2,\\[3mm]
B_3:=3\inte\Big(\frac{\Om^2}{4}|x|^2+4x_{2}^2\Big)Q\psi_{1}, \quad
B_4:=\frac{B_3}{B_1}
      +\frac{(1-4\gamma_1\gamma_2)B_1}{(2\B^*-a_1-a_2)4\gamma_1^2\gamma_2^2}.
\end{split}
\end{equation}
Here $A<0$ is given by
\begin{equation}\label{1.22:A}
A:=-\frac{\inte |x|^2Q^2dx}{a^*}<0,
\end{equation}
and $\gam_1,\gam_2\in(0,1)$ are as in \eqref{1.15:gam}.
Furthermore, we denote $\rho_{j\B}>0$ by
\begin{equation}\label{1.23:rhojB}
\ds\rho_{j\B}:=\sqrt{\frac{a^*(\B-a_m)}{\B^2-a_1a_2}}>0
\ \ \hbox{as} \ \ \B\nearrow\B^*,
\ \ j,m=1,2\ \ \hbox{and} \ \ j\neq m.
\end{equation}
Then we have the following refined spike profiles of minimizers.
\begin{thm}\label{thm3}
Suppose $V(x)$ satisfies \eqref{1.4:Vx}, and assume $0<\Om<\infty$. Let $(u_{1\B},u_{2\B})$ be the unique complex-valued minimizer of $e(\Om,\f)$, where $0< a_1,a_2<a^*$ and $0<\B<\B^*$. Then we have for $j=1,2$,
\begin{equation}\label{1.24:exp}
\begin{split}
\widetilde{v}_{j\B}(x):&=\sqrt{a^*}\alpha_\B u_{j\B}(\alpha_\B x+x_\B)
e^{-i\big(\frac{\alpha_\B\Om}{2}x\cdot x_\B^\perp-\widetilde{\theta}_{j\B}\big)}\\
&=\rho_{j\B}\Big(Q(x)+\alpha_\B^2C_{1}(x)+\alpha_\B^4C_{2}(x)\Big)+o(\alpha_\B^4)\ \
\hbox{in}\ \ L^\infty(\R^2,\C) \ \
\hbox{as}\ \ \B\nearrow\B^*,
\end{split}
\end{equation}
where $\alpha_\B>0$ is given by \eqref{1.16:alp}, $x_\B\in\R^2$ is the unique maximal point of $|u_{1\B}|^2+|u_{2\B}|^2$ and satisfies \eqref{1.18:xB}, $\widetilde{\theta}_{j\B}\in [0,2\pi)$ is a properly chosen constant, and $\rho_{j\B}>0$ is as in \eqref{1.23:rhojB}.
Here the functions $C_1(x)$ and $C_2(x)$ are defined by \eqref{1.20:Cx}, and $B_1,B_2,B_3,B_4$ are as in \eqref{1.21:B1234}.
\end{thm}


We remark that the second term $\alpha_\B^2C_{1}(x)$ of the expansion is bigger than those in \cite{GLW,GGLL} due to the degeneracy of the ring-shaped potentials.
Theorem \ref{thm3} provides a more refined characterization of minimizers for $e(\Om,\f)$ as $\B\nearrow\B^*$. Its proof needs the refined expansions of $x_\B$, $v_{j\B}(x)$ and $\eps_\B$, where $v_{j\B}(x)$ and $\eps_\B$ are as in \eqref{2.9:vjB} and \eqref{2.7:eps} below, respectively, see Lemmas \ref{lem4.1}, \ref{lem4.2} and \ref{lem4.3} for details. Moreover, we shall use the Pohozaev identity about the magnetic translation
to estimate the global maximal point $x_\B$ of $|u_{1\B}|^2+|u_{2\B}|^2$. It should point out that our method of deriving the expansion  is more efficient than \cite{GLW,GGLL}.


This paper is organized as follows. In Section 2, we shall complete the proof of Theorem \ref{thm1} on the limiting behavior of minimizers for $e(\Om,\f)$. Applying the limiting behavior of minimizers for $e(\Om,\f)$ in Theorem \ref{thm1}, we shall complete the proof of Theorem \ref{thm2} on the axial symmetry of minimizers for $e(\Om,\f)$ in Section 3.
Section 4 is devoted to proving Theorem \ref{thm3} on the refined spike profiles of minimizers for $e(\Om,\f)$ as $\B\nearrow\B^*$.

\section{Limiting Behavior of Minimizers}
In this section, we shall address the proof of Theorem \ref{thm1} on the limiting behavior of minimizers for $e(\Om,\f)$ as $\B\nearrow\B^*$, where $0<\Om<\infty$ and $0< a_1,a_2<a^*$ are fixed. Towards this purpose,
we first recall from \cite [Lemma 8.1.2] {C} and \cite [Proposition 4.1]{GNN} that the unique positive solution $Q$ of the equation \eqref{1.10:Q} satisfies
\begin{equation}\label{2.1:iden}
\inte Q^2dx=\frac{1}{2}\inte Q^4dx=\inte |\nabla Q|^2dx,
\end{equation}
and
\begin{equation} \label{2.2:dec}
Q(x) \, , \ |\nabla Q(x)| = O(|x|^{-\frac{1}{2}}e^{-|x|})
\ \ \text{as} \  \ |x|\to \infty.
\end{equation}
Recall also from \cite{Cing,Cing2} that $Q$ is non-degenerate in the sense that its linearized operator
\begin{equation}\label{2.3:L}
\mathcal{L}:=-\Delta+1-Q^2-2Re(Q\cdot)Q\ \ \hbox{in}\ \ L^2(\R^2,\C)
\end{equation}
satisfies
\begin{equation}\label{2.4:Lnon}
ker \mathcal{L}=span\l\{iQ,\, \frac{\p Q}{\p x_1},\, \frac{\p Q}{\p x_2}\r\}.
\end{equation}
We now introduce the following system
\begin{equation}\label{2.5:Q12Beq}
 \left\{\begin{array}{lll}
\ds-\Delta u_{1}+u_{1}-\frac{a_1}{a^*}u_{1}^3
                       -\frac{\B}{a^*}u_{2}^2u_{1}=0
\quad \hbox{in}\  \ \R^2,\\[3mm]
\ds-\Delta u_{2}+u_{2}-\frac{a_2}{a^*}u_{2}^3
                       -\frac{\B}{a^*}u_{1}^2u_{2}=0
\quad \hbox{in}\  \ \R^2,
\end{array}\right.
(u_1,u_2)\in H^1(\R^2,\R)\times H^1(\R^2,\R),
\end{equation}
where $a_1,a_2>0$ and $\B>\max\{a_1,a_2\}$.
It then follows from \cite [Theorem 1.3] {WY} that \eqref{2.5:Q12Beq} admits a unique positive solution $(Q_{1\B},Q_{2\B})$ (i.e., $Q_{1\B}>0$ and $Q_{2\B}>0$) satisfying
\begin{equation}\label{2.6:Q12B}
(Q_{1\B},Q_{2\B})=(\rho_{1\B}Q,\rho_{2\B}Q),
\end{equation}
where $\rho_{j\B}>0$ is given by \eqref{1.23:rhojB}, and $Q$ is the unique positive solution of the equation \eqref{1.10:Q}.

Let $(u_{1\B},u_{2\B})$ be a minimizer of $e(\Om,\f)$, where $0<\Om<\infty$, $0< a_1,a_2<a^*$ and $0<\B<\B^*$.
Recall from \eqref{1.14:mu-infty} that the Lagrange multiplier $\mu_\B$ in \eqref{1.13:mu} satisfies $\mu_\B\to-\infty$ as $\B\nearrow\B^*$.
Denote
\begin{equation}\label{2.7:eps}
\eps_\B:=\sqrt{\frac{1}{-\mu_\B}}>0\ \ \hbox{as}\ \ \B\nearrow\B^*,
\end{equation}
so that
\begin{equation}\label{2.8:eps0}
\eps_\B\to0\ \ \hbox{as}\ \ \B\nearrow\B^*.
\end{equation}
Define
\begin{equation}\label{2.9:vjB}
v_{j\B}(x):=\sqrt{a^*}\eps_\B u_{j\B}\big(\eps_\B x+x_\B\big)
    e^{-i\big(\frac{\eps_\B\Om}{2}x\cdot x_\B^\perp-\theta_{j\B}\big)},
\,\ j=1,2,
\end{equation}
where $x_\B\in\R^2$ is a global maximal point of $|u_{1\B}|^2+|u_{2\B}|^2$, and the constant phase $\theta_{j\B}\in [0,2\pi)$ is chosen such that
\begin{equation}\label{2.10:theta}
\big\|v_{j\B}-\sqrt{\gamma_j}Q\big\|_{L^2(\R^2)}
=\min\limits_{\theta\in[0,2\pi)}\big\|
\sqrt{a^*}\eps_\B u_{j\B}(\eps_\B x+x_\B)
e^{-i\big(\frac{\eps_\B\Om}{2}x\cdot x_\B^\perp-\theta\big)}
-\sqrt{\gamma_j}Q\big\|_{L^2(\R^2)},
\end{equation}
where $\gam_j\in(0,1)$ is given by \eqref{1.15:gam}.
It then follows from \eqref{2.10:theta} that
\begin{equation}\label{2.11:ortho}
Re\Big(\inte v_{j\B}(iQ)\Big)=0, \  \  j=1,2.
\end{equation}
Similar to the arguments of proving \cite [Lemmas 3.4 and 3.5] {GGLL}, one can derive from \eqref{2.8:eps0}--\eqref{2.11:ortho} that
\begin{equation}\label{2.12:lim}
v_{j\B}(x)\to\sqrt{\gamma_j}Q(x)
\ \ \hbox{strongly in} \ \ H^1(\R^2,\C)\cap L^\infty(\R^2,\C)
\ \ \hbox{as}\ \ \B\nearrow\B^*, \ \ j=1,2,
\end{equation}
and the global maximal point $x_\B\in\R^2$ of $|u_{1\B}|^2+|u_{2\B}|^2$ is unique and satisfies
\begin{equation}\label{2.13:xB}
x_\B\to P_0\ \ \hbox{as}\ \ \B\nearrow\B^*
\end{equation}
for some $P_0\in\R^2$ satisfying $|P_0|=1$. Note that the system \eqref{1.1:GPs} is invariant under the rotational transformation. Without loss of generality, we may assume that $x_{\B}$ lies on the $x_2$-axis, i.e., $x_{\B}=(0,p_{\B})$.
Following \eqref{1.12:ELs}, one can deduce that $\big(v_{1\B}(x),v_{2\B}(x)\big)$ defined by \eqref{2.9:vjB} satisfies the following system
\begin{equation}\label{2.14:GPs}
\left\{\begin{array}{lll}
\ds-\Big(\nabla-i\frac{\eps_\B^2\Om}{2}x^\perp\Big)^2v_{1\B}
+\Big[\eps_\B^2\big(|\eps_\B x+x_\B|^2-1\big)^2\\[3mm]
+1-\ds \frac{a_1}{a^*}|v_{1\B}|^2-\frac{\B}{a^*}|v_{2\B}|^2
\Big]v_{1\B}=0\quad \hbox{in}\  \ \R^2,\\[3mm]
\ds-\Big(\nabla-i\frac{\eps_\B^2\Om}{2}x^\perp\Big)^2v_{2\B}
+\Big[\eps_\B^2\big(|\eps_\B x+x_\B|^2-1\big)^2\\[3mm]
+1-\ds\frac{a_2}{a^*}|v_{2\B}|^2-\frac{\B}{a^*}|v_{1\B}|^2
\Big]v_{2\B}=0\quad \hbox{in}\  \ \R^2.
\end{array}\right.
\end{equation}
By the comparison principle, one can derive from \eqref{2.12:lim}--\eqref{2.14:GPs} that
\begin{equation}\label{2.15:dec}
|v_{j\B}(x)|\le Ce^{-\frac{2}{3}|x|}, \ \ |\nabla v_{j\B}(x)|\le Ce^{-\frac{1}{2}|x|}
\ \ \hbox{in}\ \ \R^2 \ \ \hbox{as}\ \ \B\nearrow\B^*, \ \ j=1,2,
\end{equation}
where the constant $C>0$ is independent of $0<\B<\B^*$.

We now rewrite $v_{j\B}(x)$ as
\begin{equation}\label{2.16:exp}
v_{j\B}(x):=\sqrt{a^*}\eps_\B u_{j\B}\big(\eps_\B x+x_\B\big)
    e^{-i\big(\frac{\eps_\B\Om}{2}x\cdot x_\B^\perp-\theta_{j\B}\big)}
    =\rho_{j\B}Q(x)+\hat{u}_{j\B}(x),\ \ j=1,2,
\end{equation}
where $\rho_{j\B}>0$ is as in \eqref{1.23:rhojB}.
Note that $\rho_{j\B}\to\sqrt{\gam_j}$ as $\B\nearrow\B^*$, where $\gam_j\in(0,1)$ is defined by \eqref{1.15:gam}. We then obtain from \eqref{2.12:lim} that
\begin{equation}\label{2.17:low}
\|\hat{u}_{j\B}(x)\|_{L^\infty(\R^2)}\to0 \ \ \hbox{as} \ \ \B\nearrow\B^*,\ \ j=1,2.
\end{equation}
Following \eqref{2.5:Q12Beq}, \eqref{2.6:Q12B} and \eqref{2.14:GPs}, we derive from  \eqref{2.16:exp} that $(\hat{u}_{1\B},\hat{u}_{2\B})$ satisfies the following system
\begin{equation}\label{2.18:GPs}
 \left\{\begin{array}{lll}
\mathcal{L}_{1\B}(\hat{u}_{1\B},\hat{u}_{2\B})=F_{1\B}(x)
\  \ \hbox{in}\  \ \R^2,\\[3mm]
\mathcal{L}_{2\B}(\hat{u}_{2\B},\hat{u}_{1\B})=F_{2\B}(x)
\  \ \hbox{in}\  \ \R^2,
\end{array}\right.
\end{equation}
where the operators $\mathcal{L}_{1\B}$ and $\mathcal{L}_{2\B}$ are given by
\begin{equation}\label{2.19:linop}
 \left\{\begin{array}{lll}
 \begin{split}
\mathcal{L}_{1\B}(\phi_1,\phi_2)
:=&-\Delta \phi_1+\phi_1
   -\frac{a_1}{a^*}|Q_{1\B}|^2\phi_1-\frac{2a_1}{a^*}Re(Q_{1\B}\phi_1)Q_{1\B}\\[3mm]
  &-\frac{\B}{a^*}|Q_{2\B}|^2\phi_1-\frac{2\B}{a^*}Re(Q_{2\B}\phi_2)Q_{1\B}
\  \ \hbox{in}\  \ \R^2,\\[3mm]
\mathcal{L}_{2\B}(\phi_2,\phi_1)
:=&-\Delta \phi_2+\phi_2
   -\frac{a_2}{a^*}|Q_{2\B}|^2\phi_2-\frac{2a_2}{a^*}Re(Q_{2\B}\phi_2)Q_{2\B}\\[3mm]
  &-\frac{\B}{a^*}|Q_{1\B}|^2\phi_2-\frac{2\B}{a^*}Re(Q_{1\B}\phi_1)Q_{2\B}
\  \ \hbox{in}\  \ \R^2,
 \end{split}
\end{array}\right.
\end{equation}
and $F_{j\B}(x)$ is defined by
\begin{equation}\label{2.20:LF}
\begin{split}
F_{j\B}(x)
:=&-\Big(\frac{\eps_\B^4\Omega^2}{4}|x|^2
                   +\eps_\B^2\big(|\eps_\B x+x_\B|^2-1\big)^2\Big)\rho_{j\B}Q
             -i\,\eps_\B^2\,\Omega \big(x^{\perp}\cdot\nabla\hat{u}_{j\B}\big)\\[3mm]
 \quad&-\Big(\frac{\eps_\B^4\Omega^2}{4}|x|^2
        +\eps_\B^2\big(|\eps_\B x+x_\B|^2-1\big)^2\Big)\hat{u}_{j\B}\\[3mm]
 \quad&+\frac{a_j}{a^*}|\hat{u}_{j\B}|^2\hat{u}_{j\B}
         +\frac{2a_j}{a^*}Re\big(\rho_{j\B}Q\hat{u}_{j\B}\big)\hat{u}_{j\B}
         +\frac{a_j}{a^*}|\hat{u}_{j\B}|^2\rho_{j\B}Q\\[3mm]
 \quad&+\frac{\B}{a^*}|\hat{u}_{m\B}|^2\hat{u}_{j\B}
         +\frac{2\B}{a^*}Re\big(\rho_{m\B}Q\hat{u}_{m\B}\big)\hat{u}_{j\B}
         +\frac{\B}{a^*}|\hat{u}_{m\B}|^2\rho_{j\B}Q,
\end{split}
\end{equation}
where $j,m=1,2$ and $j\neq m$.

\subsection{Refined estimates of  $(v_{1\B},v_{2\B})$}
In this subsection, we shall derive the refined estimates of $(v_{1\B},v_{2\B})$ defined in \eqref{2.9:vjB}.
We first give the estimates of $\hat{u}_{j\B}(x)$ given by \eqref{2.16:exp} as $\B\nearrow\B^*$ for $j=1,2$.
\begin{lem}\label{lem2.1}
Suppose $V(x)$ satisfies \eqref{1.4:Vx}, and assume $0<\Om<\infty$ and
$0< a_1,a_2<a^*$ are fixed. Then the function $\hat{u}_{j\B}(x)$ in \eqref{2.16:exp} satisfies
\begin{equation}\label{2.21:lowdec}
|\hat{u}_{j\B}(x)|\leq C\eps_\B^2e^{-\frac{1}{4}|x|}, \ \
|\nabla\hat{u}_{j\B}(x)|\leq C\eps_\B^2e^{-\frac{1}{8}|x|}
\ \ \hbox{uniformly in}\ \ \R^2
\ \ \hbox{as}\ \ \B\nearrow\B^*, \ \ j=1,2,
\end{equation}
where the constant $C>0$ is independent of $0<\B<\B^*$.
\end{lem}
\noindent{\bf Proof.} We first claim that
\begin{equation}\label{2.22}
\|\hat{u}_{j\B}(x)\|_{L^\infty(\R^2)}\leq C\eps_\B^2
\ \ \hbox{uniformly in}\ \ \R^2
\ \ \hbox{as}\ \ \B\nearrow\B^*, \ \ j=1,2,
\end{equation}
where the constant $C>0$ is independent of $0<\B<\B^*$.
By contradiction, we assume that
\begin{equation}\label{2.23}
\lim\limits_{\B\nearrow\B^*}\frac{\|\hat{u}_{j\B}\|_{L^\infty(\R^2)}}{\eps_\B^2}=\infty, \ \ j=1,2.
\end{equation}
Define
\begin{equation}\label{2.24:UjB}
U_{j\B}(x):=\frac{\hat{u}_{j\B}(x)}
{\max\{\|\hat{u}_{1\B}\|_{L^\infty(\R^2)},\|\hat{u}_{2\B}\|_{L^\infty(\R^2)}\}}, \ \ j=1,2.
\end{equation}
We then obtain from \eqref{2.18:GPs} and \eqref{2.24:UjB} that $(U_{1\B},U_{2\B})$ satisfies
\begin{equation}\label{2.25:GPs}
 \left\{\begin{array}{lll}
\mathcal{L}_{1\B}(U_{1\B},U_{2\B})
=\ds\frac{F_{1\B}(x)}
{\max\{\|\hat{u}_{1\B}\|_{L^\infty(\R^2)},\|\hat{u}_{2\B}\|_{L^\infty(\R^2)}\}}
\  \ \hbox{in}\  \ \R^2,\\[3mm]
\mathcal{L}_{2\B}(U_{2\B},U_{1\B})
=\ds\frac{F_{2\B}(x)}
{\max\{\|\hat{u}_{1\B}\|_{L^\infty(\R^2)},\|\hat{u}_{2\B}\|_{L^\infty(\R^2)}\}}
\  \ \hbox{in}\  \ \R^2,
\end{array}\right.
\end{equation}
where the operators $\mathcal{L}_{1\B}$ and $\mathcal{L}_{2\B}$ are as in \eqref{2.19:linop}, and $F_{j\B}(x)$ is given by \eqref{2.20:LF} for $j=1,2$.
Using the comparison principle, one can derive from \eqref{2.2:dec}, \eqref{2.18:GPs}--\eqref{2.20:LF} that
\begin{equation}\label{2.26:lowdec}
|\hat{u}_{j\B}(x)|\le C\delta_\B e^{-\frac{2}{3}|x|}, \ \ |\nabla \hat{u}_{j\B}(x)|\le C\delta_\B e^{-\frac{1}{2}|x|}
\ \ \hbox{in}\ \ \R^2 \ \ \hbox{as}\ \ \B\nearrow\B^*, \ \ j=1,2,
\end{equation}
where $\delta_\B$ satisfies $\delta_\B=o(1)$ as $\B\nearrow\B^*$.
By \eqref{2.2:dec} and \eqref{2.26:lowdec}, we deduce from \eqref{2.20:LF} that for $j=1,2$,
\begin{equation}\label{2.27}
\begin{split}
&\quad
\frac{|F_{j\B}(x)|}
{\max\{\|\hat{u}_{1\B}\|_{L^\infty(\R^2)},\|\hat{u}_{2\B}\|_{L^\infty(\R^2)}\}}\\
&\le \frac{C\eps_\B^2}
{\max\{\|\hat{u}_{1\B}\|_{L^\infty(\R^2)},\|\hat{u}_{2\B}\|_{L^\infty(\R^2)}\}}\\
&\quad+C
\max\{\|\hat{u}_{1\B}\|_{L^\infty(\R^2)},\|\hat{u}_{2\B}\|_{L^\infty(\R^2)}\}
e^{-\frac{2}{3}|x|}
\ \ \hbox{uniformly in}\ \ \R^2
\ \ \hbox{as}\ \ \B\nearrow \B^*.
\end{split}
\end{equation}
Note from \eqref{2.24:UjB} that $\|U_{j\B}(x)\|_{L^\infty(\R^2)}\le1$.
Using \eqref{2.17:low}, \eqref{2.23} and \eqref{2.27}, the elliptic regularity theory then yields from \eqref{2.25:GPs} that there exists a constant $C>0$, independent of $0<\B<\B^*$, such that $\|U_{j\B}\|_{C^{2,\alpha}_{loc}(\R^2,\C)}\le C$ for $\alpha\in(0,1)$ and $j=1,2$. Therefore, one can deduce that up to a subsequence if necessary,
\begin{equation}\label{2.28}
U_{j\B}\to U_{j0}\ \ \hbox{in} \ \ C^{2}_{loc}(\R^2,\C)
\ \ \hbox{as}\ \ \B\nearrow \B^*, \ \ j=1,2,
\end{equation}
where $(U_{10},U_{20})$ satisfies the following system
\begin{equation}\label{2.29}
 \left\{\begin{array}{lll}
 \begin{split}
\mathcal{L}_{1}(U_{10},U_{20})
:=&-\Delta U_{10}+U_{10}
   -\frac{a_1\gamma_1}{a^*}Q^2U_{10}-\frac{2a_1\gamma_1}{a^*}Re(QU_{10})Q\\[3mm]
  &-\frac{\B^*\gamma_2}{a^*}Q^2U_{10}
   -\frac{2\B^*\sqrt{\gamma_1\gamma_2}}{a^*}Re(QU_{20})Q=0
\  \ \hbox{in}\  \ \R^2,\\[3mm]
\mathcal{L}_{2}(U_{20},U_{10})
:=&-\Delta U_{20}+U_{20}
   -\frac{a_2\gamma_2}{a^*}Q^2U_{20}-\frac{2a_2\gamma_2}{a^*}Re(QU_{20})Q\\[3mm]
  &-\frac{\B^*\gamma_1}{a^*}Q^2U_{20}
   -\frac{2\B^*\sqrt{\gamma_1\gamma_2}}{a^*}Re(QU_{10})Q=0
\  \ \hbox{in}\  \ \R^2.
 \end{split}
\end{array}\right.
\end{equation}
Here $\gam_1,\gam_2\in(0,1)$ are as in \eqref{1.15:gam}. Following \cite [Lemma 2.2 and Theorem 3.1] {DW}, one can derive from \eqref{2.29} that
\begin{equation}\label{2.30}
\begin{pmatrix}
  U_{10}\\[5mm]
  U_{20}
\end{pmatrix}
=\begin{pmatrix}
ib_0^1\sqrt{\gamma_1}Q\\[5mm]
ib_0^2\sqrt{\gamma_2}Q
\end{pmatrix}
+
\sum^2_{l=1}b_l
\begin{pmatrix}
\ds\sqrt{\gamma_1}\frac{\partial Q}{\partial x_l}\\[5mm]
\ds\sqrt{\gamma_2}\frac{\partial Q}{\partial x_l}
\end{pmatrix},
\end{equation}
where $b_0^1,b_0^2,b_1,b_2\in\R$ are some constants.

We now prove that the constants $b_0^j,b_j$ in \eqref{2.30} satisfy
\begin{equation}\label{2.31}
b_0^j=b_j=0,  \ \ j=1,2.
\end{equation}
By \eqref{2.11:ortho}, we derive from \eqref{2.16:exp} and \eqref{2.24:UjB} that
$$Re\Big(\inte U_{j\B}(iQ)\Big)=0, \ \ j=1,2.$$
We then deduce from \eqref{2.28}, \eqref{2.30} and above that
$$0=Re\Big(\inte U_{j0}(iQ)\Big)=-b_0^j\sqrt{\gam_j}\inte Q^2, \ \ j=1,2,$$
which implies that $b_0^j=0$ for $j=1,2$. Note that $x=0$ is the unique global maximal point of $|v_{1\B}|^2+|v_{2\B}|^2$, we then obtain from \eqref{2.16:exp} that
$$0=\Big(\nabla\big(|v_{1\B}|^2+|v_{2\B}|^2\big)\Big)(0)
=\sum^2_{j=1}
\Big(\nabla\big(\rho_{j\B}^2Q^2+2Re(\hat{u}_{j\B}\rho_{j\B}Q)+|\hat{u}_{j\B}|^2\big)\Big)(0),$$
which implies that
\begin{equation}\label{2.32}
\sum^2_{j=1}Re(\rho_{j\B}\nabla\hat{u}_{j\B}(0))
=-\frac{\sum^2\limits_{j=1}Re(\nabla\hat{u}_{j\B}(0)\overline{\hat{u}_{j\B}}(0))}{Q(0)}.
\end{equation}
Using \eqref{2.26:lowdec} and \eqref{2.32}, we derive from \eqref{2.24:UjB} that
\begin{equation}\label{2.33}
\begin{split}
\sum^2_{j=1}Re(\rho_{j\B}\nabla U_{j\B}(0))
=-\frac{\sum^2\limits_{j=1}Re(\nabla\hat{u}_{j\B}(0)\overline{\hat{u}_{j\B}}(0))}
{Q(0)\cdot\max\{\|\hat{u}_{1\B}\|_{L^\infty(\R^2)},\|\hat{u}_{2\B}\|_{L^\infty(\R^2)}\}}
\to0
\ \ \hbox{as}\ \ \B\nearrow \B^*.
\end{split}
\end{equation}
Combining \eqref{2.28}, \eqref{2.30} with \eqref{2.33} then yields that
$$(\gamma_1+\gamma_2)
\sum^2_{l=1}b_l\frac{\p^2Q(0)}{\p x_l\p x_k}=0, \ \ k=1,2,$$
where we also have used the fact that $\rho_{j\B}\to\sqrt{\gam_j}>0$ as $\B\nearrow\B^*$.
Since $det\big(\frac{\p^2Q(0)}{\p x_l\p x_k}\big)\neq0$, we derive from above that $b_1=b_2=0$. Therefore, we conclude from above that \eqref{2.31} holds true.

Following \eqref{2.30} and \eqref{2.31}, we obtain that $U_{j0}\equiv0$ in $\R^2$, where $j=1,2$. On the other hand, using the comparison principle, we derive from \eqref{2.23}, \eqref{2.25:GPs} and \eqref{2.27} that there exists a constant $C>0$, independent of $0<\B<\B^*$, such that
\begin{equation}\label{2.34}
|U_{j\B}(x)|\le Ce^{-\frac{1}{4}|x|}
\ \ \hbox{uniformly in}\ \ \R^2
\ \ \hbox{as} \ \ \B\nearrow\B^*, \ \ j=1,2.
\end{equation}
Taking $x_\B\in\R^2$ such that
\begin{equation}\label{2.35}
\max\{U_{1\B}(x_\B),U_{2\B}(x_\B)\}=1.
\end{equation}
By \eqref{2.34}, we deduce that $|x_\B|\le C$ uniformly in $\B$. Using \eqref{2.28}, we then obtain that $$
\max\{\|U_{10}\|_{L^\infty(\R^2)},\|U_{20}\|_{L^\infty(\R^2)}\}\ge1,$$
which contradicts with $U_{j0}\equiv0$ in $\R^2$, and hence \eqref{2.22} holds true.

Using \eqref{2.2:dec}, \eqref{2.15:dec} and \eqref{2.22}, we derive from \eqref{2.20:LF} that
\begin{equation}\label{2.36:GPs}
|F_{j\B}(x)|\le C\eps_\B^2e^{-\frac{1}{2}|x|}
\ \ \hbox{uniformly in}\ \ \R^2
\  \ \hbox{as}\ \ \B\nearrow\B^*, \ \ j=1,2.
\end{equation}
By the comparison principle, we deduce from \eqref{2.18:GPs}, \eqref{2.22} and \eqref{2.36:GPs} that
\begin{equation}\label{2.37}
|\hat{u}_{j\B}(x)|\le C\eps_\B^2e^{-\frac{1}{4}|x|}
\ \ \hbox{uniformly in}\ \ \R^2
\ \ \hbox{as} \ \ \B\nearrow\B^*, \ \ j=1,2.
\end{equation}
Moreover, using the comparison principle again, one can derive from \eqref{2.18:GPs} and \eqref{2.37} that
\begin{equation}\label{2.38}
|\nabla\hat{u}_{j\B}(x)|\le C\eps_\B^2e^{-\frac{1}{8}|x|}
\ \ \hbox{uniformly in}\ \ \R^2
\ \ \hbox{as} \ \ \B\nearrow\B^*, \ \ j=1,2.
\end{equation}
It then follows from \eqref{2.37} and \eqref{2.38} that \eqref{2.21:lowdec} holds true. Lemma \ref{lem2.1} is thus proved.
\qed

\vskip 0.1truein

We next establish the expansion of $x_\B$ in terms of $\eps_\B$, where $x_\B$ is the unique maximal point of $|u_{1\B}|^2+|u_{2\B}|^2$.
\begin{lem}\label{lem2.2}
Suppose $V(x)$ satisfies \eqref{1.4:Vx}, and assume $0<\Om<\infty$ and
$0< a_1,a_2<a^*$ are fixed.
Let $(u_{1\B},u_{2\B})$ be a minimizer of $e(\Om,\f)$ and $x_\B=(0,p_\B)$ be the unique maximal point of $|u_{1\B}|^2+|u_{2\B}|^2$. Then we have
\begin{equation}\label{2.39:pB}
p_{\B}-1=A\eps_\B^2+o(\eps_\B^2)
\ \ \hbox{as}\ \ \B\nearrow\B^*,
\end{equation}
where $A=-\frac{\inte |x|^2Q^2dx}{a^*}<0$ is as in \eqref{1.22:A}.
\end{lem}
\noindent{\bf Proof.}
Similar to the argument of proving \cite [Lemma A.1] {GGLL}, one can derive from \eqref{2.14:GPs} and \eqref{2.15:dec} that when $\B$ close to $\B^*$ enough,
\begin{equation}\label{2.40}
\inte\frac{\p \big(|\eps_\B x+x_\B|^2-1\big)^2}
{\p x_{2}}\big(|v_{1\B}|^2+|v_{2\B}|^2\big)=0.
\end{equation}

On the other hand, we derive from \eqref{2.16:exp} and Lemma \ref{lem2.1} that
\begin{equation}\label{2.41}
\begin{aligned}
&\quad\sum^2_{j=1}\inte\frac{\p \big(|\eps_\B x+x_\B|^2-1\big)^2}
{\p x_{2}}|v_{j\B}|^2dx\\
&=\sum^2_{j=1}\inte4\eps_\B(|\eps_\B x+x_\B|^2-1)(\eps_\B x_2+p_\B)
                         \big|\rho_{j\B}Q+\hat{u}_{j\B}\big|^2dx\\
&=\sum^2_{j=1}\inte4\eps_\B(\eps_\B^2|x|^2+2\eps_\B x_2p_\B+p_\B^2-1)
                           (\eps_\B x_2+p_\B)
             \l(\rho_{j\B}^2Q^2+2Re(\rho_{j\B}Q\hat{u}_{j\B})+|\hat{u}_{j\B}|^2\r)dx\\
&=\sum^2_{j=1}\rho_{j\B}^2\Big[8\eps_\B^3\inte|x|^2Q^2dx
                              +8\eps_\B(p_{\B}-1)\inte Q^2dx\Big]
  +o(\eps_\B|p_{\B}-1|)+o(\eps_\B^3)
\ \ \hbox{as}\ \ \B\nearrow\B^*.
\end{aligned}
\end{equation}
Combining \eqref{2.40} with \eqref{2.41} then yields that \eqref{2.39:pB} holds true, and Lemma \ref{lem2.2} is therefore proved.\qed

\vskip 0.1truein

Based on Lemmas \ref{lem2.1} and \ref{lem2.2},
we now establish the expansion of $\big(v_{1\B}(x),v_{2\B}(x)\big)$ in terms of $\eps_\B$.
\begin{lem}\label{lem2.3}
Under the assumptions of Lemma \ref{lem2.2}, and let $v_{j\B}(x)$ be defined by \eqref{2.9:vjB}, where $j=1,2$. Then we have
\begin{equation}\label{2.42}
\begin{aligned}
v_{j\B}(x)
=\rho_{j\B}Q(x)+\rho_{j\B}\eps_\B^4\psi_{1}(x)+o(\eps_\B^4)\ \ \hbox{in}\ \ \R^2
\ \ \hbox{as} \ \ \B\nearrow\B^*, \ \ j=1,2,
\end{aligned}
\end{equation}
where $\rho_{j\B}>0$ is as in \eqref{1.23:rhojB}, and $\psi_{1}(x)\in C^2(\R^2,\C)\cap L^\infty(\R^2,\C)$ solves uniquely
\begin{equation}\label{2.43:psi1}
\mathcal{L}\psi_{1}(x)=-\Big(\frac{\Om^2}{4}|x|^2+4x_{2}^2\Big)Q(x)\ \ \mbox{in}\ \ \R^2,
\ \ \nabla Re(\psi_{1}(0))=0,\ Re\inte \psi_1(iQ)=0,
\end{equation}
where the operator $\mathcal{L}$ is defined by \eqref{2.3:L}.
\end{lem}
\noindent{\bf Proof.}
Using Lemmas \ref{lem2.1} and \ref{lem2.2}, we derive from \eqref{2.20:LF} that
\begin{equation}\label{2.44:GPs}
|F_{j\B}(x)|\le C\eps_\B^4e^{-\frac{1}{8}|x|}
\ \ \hbox{uniformly in}\ \ \R^2
\  \ \hbox{as}\ \ \B\nearrow\B^*, \ \ j=1,2.
\end{equation}
Similar to the argument of proving Lemma \ref{lem2.1}, by the comparison principle, one can deduce from \eqref{2.18:GPs} and \eqref{2.44:GPs} that
\begin{equation}\label{2.45}
|\hat{u}_{j\B}(x)|\leq C\eps_\B^4e^{-\frac{1}{8}|x|}, \ \
|\nabla\hat{u}_{j\B}(x)|\leq C\eps_\B^4e^{-\frac{1}{16}|x|}
\ \ \hbox{uniformly in}\ \ \R^2
\ \ \hbox{as}\ \ \B\nearrow\B^*, \ \ j=1,2,
\end{equation}
where the constant $C>0$ is independent of $0<\B<\B^*$.

Denote
\begin{equation}\label{2.46}
\mathcal{U}_{j\B}(x):=\hat{u}_{j\B}(x)-\rho_{j\B}\eps_\B^4\psi_{1}(x),\ \ j=1,2,
\end{equation}
where $\rho_{j\B}>0$ is defined by \eqref{1.23:rhojB}, $\hat{u}_{j\B}(x)$ is as in \eqref{2.16:exp}, and
$\psi_{1}(x)\in C^2(\R^2,\C)\cap L^\infty(\R^2,\C)$ is a solution of
\eqref{2.43:psi1}.
It then follows from \eqref{2.4:Lnon} and \eqref{2.43:psi1} that $\psi_{1}(x)$ is unique.
Following \eqref{2.18:GPs} and \eqref{2.43:psi1}, we obtain that $(\mathcal{U}_{1\B},\mathcal{U}_{2\B})$ defined by \eqref{2.46} satisfies
\begin{equation}\label{2.47}
\left\{\begin{array}{lll}
\mathcal{L}_{1\B}(\mathcal{U}_{1\B},\mathcal{U}_{2\B})
 =F_{1\B}(x)+\rho_{1\B}\eps_\B^4\big(\frac{\Om^2}{4}|x|^2+4x_{2}^2\big)Q
:=\mathcal{F}_{1\B}(x)
\quad \mbox{in}\,\  \R^2,\\[3mm]
\mathcal{L}_{2\B}(\mathcal{U}_{2\B},\mathcal{U}_{1\B})
 =F_{2\B}(x)+\rho_{2\B}\eps_\B^4\big(\frac{\Om^2}{4}|x|^2+4x_{2}^2\big)Q
:=\mathcal{F}_{2\B}(x)
\quad \mbox{in}\,\  \R^2.
\end{array}\right.
\end{equation}
By \eqref{2.2:dec} and Lemma \ref{lem2.2}, we deduce from \eqref{2.45} that the terms $\mathcal{F}_{1\B}(x)$ and $\mathcal{F}_{2\B}(x)$ satisfy
\begin{equation}\label{2.48}
\left\{\begin{array}{lll}
\begin{split}
 \ds|\mathcal{F}_{1\B}(x)|:&=\Big|F_{1\B}(x)
                  +\rho_{1\B}\eps_\B^4\Big(\frac{\Om^2}{4}|x|^2+4x_{2}^2\Big)Q\Big|
 \le C_{1\B}\eps_\B^4e^{-\frac{1}{16}|x|}\\[3mm]
 \ds|\mathcal{F}_{2\B}(x)|:&=\Big|F_{2\B}(x)
                  +\rho_{2\B}\eps_\B^4\Big(\frac{\Om^2}{4}|x|^2+4x_{2}^2\Big)Q\Big|
 \le C_{2\B}\eps_\B^4e^{-\frac{1}{16}|x|}
\end{split}
\end{array}\right.
\ \ \hbox{uniformly in}\ \ \R^2\ \ \hbox{as}\ \  \B\nearrow \B^*,
\end{equation}
where $C_{j\B}>0$ satisfies $C_{j\B}=o(1)$ as $\B\nearrow\B^*$ for $j=1,2$.
Similar to \eqref{2.22} in Lemma \ref{lem2.1}, one can derive from \eqref{2.45}, \eqref{2.47} and \eqref{2.48} that
\begin{equation}\label{2.49}
\big| \mathcal{U}_{j\B}(x)\big|
\leq C_{j\B} \eps_\B^4\quad \hbox{uniformly in}\ \ \R^2
\ \ \hbox{as}\, \ \B\nearrow \B^*,\ \  j=1,2,
\end{equation}
where $C_{j\B}>0$ satisfies $C_{j\B}=o(1)$ as $\B\nearrow\B^*$. Combining \eqref{2.46} with \eqref{2.49} yields that
\begin{equation}\label{2.50}
\hat{u}_{j\B}(x)=\rho_{j\B}\eps_\B^4\psi_{1}(x)+o(\eps_\B^{4})
\ \ \hbox{in} \ \ \R^2
\ \ \hbox{as}\ \  \B\nearrow \B^*,\ \ j=1,2.
\end{equation}
It then follows from \eqref{2.16:exp} and \eqref{2.50} that \eqref{2.42} holds true, and the proof of Lemma \ref{lem2.3} is therefore complete.
\qed

\subsection{Proof of Theorem \ref{thm1}}
The main purpose of this subsection is to complete the proof of Theorem \ref{thm1}. We first give the blow up rate of minimizers for $e(\Om,\f)$ in the following lemma.

\begin{lem}\label{lem2.4}
Under the assumptions of Lemma \ref{lem2.2}, the term $\eps_\B$ in \eqref{2.7:eps} satisfies
\begin{equation}\label{2.51:epal}
  \eps_\B=\alpha_\B+o(\alpha_\B)\ \ \hbox{as}\ \ \B\nearrow\B^*,
\end{equation}
where $\alpha_\B>0$ is as in \eqref{1.16:alp}.
\end{lem}
\noindent{\bf Proof.}
Multiplying the first equation of \eqref{2.14:GPs} by
$x\cdot\nabla\ovl{v_{1\B}}$ and the second equation of \eqref{2.14:GPs} by
$x\cdot\nabla\ovl{v_{2\B}}$, integrating over $\R^2$ and taking theirs real parts, we derive from \eqref{2.15:dec} that $(v_{1\B},v_{2\B})$ defined by \eqref{2.9:vjB} satisfies the following Pohozaev identity
\begin{equation}\label{2.52}
\begin{split}
&\quad\sum^2_{j=1}\inte\Big(\frac{a_j}{2a^*}|v_{j\B}|^4-|v_{j\B}|^2\Big)
    +\frac{\B}{a^*}\inte|v_{1\B}|^2|v_{2\B}|^2\\
&=\sum^2_{j=1}\inte\Big(\frac{\eps_\B^4\Omega^2}{2}|x|^2
                +\eps_\B^2\big(|\eps_\B x+x_\B|^2-1\big)^2\\
  &\qquad\qquad\quad+2\eps_\B^3(|\eps_\B x+x_\B|^2-1)
                (\eps_\B|x|^2+x\cdot x_\B)\Big)|v_{j\B}|^2\\
  &\quad-\sum^2_{j=1}\inte Re\l[i\eps_\B^2\Om(x^\perp\cdot\nabla v_{j\B})(x\cdot\nabla\ovl{v_{j\B}})\r].
\end{split}
\end{equation}

We shall first estimate the left hand side of \eqref{2.52}. Using \eqref{1.23:rhojB} and \eqref{2.1:iden}, we deduce from \eqref{2.16:exp} and Lemma \ref{lem2.1} that
\begin{equation}\label{2.53}
\begin{split}
&\quad L.\ H.\ S.\ \ \hbox{of}\ \ \eqref{2.52}\\
&=\frac{(a^*)^2(2\B-a_1-a_2)}{\B^2-a_1a_2}-a^*\\
&\quad+\sum^2_{j=1}\inte\frac{a_j}{2a^*}
\Big(2\rho_{j\B}^2Q^2|\hat{u}_{j\B}|^2
   +4Re(\rho_{j\B}^3Q^3\hat{u}_{j\B})
   +4\big(Re(\rho_{j\B}Q\hat{u}_{j\B})\big)^2\Big)\\
&\quad+\frac{\B}{a^*}\inte\Big(
    \rho_{1\B}^2Q^2|\hat{u}_{2\B}|^2+\rho_{2\B}^2Q^2|\hat{u}_{1\B}|^2
   +2Re(\rho_{1\B}^2\rho_{2\B}Q^3\hat{u}_{2\B})+2Re(\rho_{1\B}\rho_{2\B}^2Q^3\hat{u}_{1\B})\\
   &\qquad\qquad\qquad+4Re(\rho_{1\B}\rho_{2\B}Q^2\hat{u}_{1\B}\hat{u}_{2\B})\Big)
+o(\eps_\B^{4})
\ \ \hbox{as}\ \ \B\nearrow\B^*,
\end{split}
\end{equation}
where we also have used the fact that
$\inte\big(|v_{1\B}|^2+|v_{2\B}|^2\big)=a^*\inte\big(|u_{1\B}|^2+|u_{2\B}|^2\big)
=a^*$.
Furthermore, we obtain from \eqref{2.5:Q12Beq}, \eqref{2.6:Q12B} and \eqref{2.18:GPs} that
\begin{equation*}
\begin{split}
&\quad\sum^2_{j=1}\inte\frac{2a_j}{a^*}Re(\rho_{j\B}^3Q^3\hat{u}_{j\B})
   +\frac{2\B}{a^*}\inte\Big(Re(\rho_{1\B}^2\rho_{2\B}Q^3\hat{u}_{2\B})
                            +Re(\rho_{1\B}\rho_{2\B}^2Q^3\hat{u}_{1\B})\Big)\\
&=2\sum^2_{j=1}\inte Re\big((-\Delta(\rho_{j\B}Q)+\rho_{j\B}Q)\hat{u}_{j\B}\big)
=2\sum^2_{j=1}\inte Re\big((-\Delta \hat{u}_{j\B}
                             +\hat{u}_{j\B})\rho_{j\B}Q\big)\\
&=3\Big\{\sum^2_{j=1}\inte\frac{2a_j}{a^*}Re(\rho_{j\B}^3Q^3\hat{u}_{j\B})
   +\frac{2\B}{a^*}\inte\Big(Re(\rho_{1\B}^2\rho_{2\B}Q^3\hat{u}_{2\B})
                            +Re(\rho_{1\B}\rho_{2\B}^2Q^3\hat{u}_{1\B})\Big)\Big\}\\
&\quad+2\sum^2_{j=1}\inte Re(F_{j\B}\rho_{j\B}Q),
\end{split}
\end{equation*}
which implies that
\begin{equation}\label{2.54}
\begin{split}
&\quad\sum^2_{j=1}\inte\frac{2a_j}{a^*}Re(\rho_{j\B}^3Q^3\hat{u}_{j\B})
   +\frac{2\B}{a^*}\inte\Big(Re(\rho_{1\B}^2\rho_{2\B}Q^3\hat{u}_{2\B})
                            +Re(\rho_{1\B}\rho_{2\B}^2Q^3\hat{u}_{1\B})\Big)\\
&=-\sum^2_{j=1}\inte Re(F_{j\B}\rho_{j\B}Q).
\end{split}
\end{equation}
We then derive from \eqref{2.20:LF} and \eqref{2.54} that
\begin{equation}\label{2.55}
\begin{split}
&\quad\sum^2_{j=1}\inte\frac{2a_j}{a^*}Re(\rho_{j\B}^3Q^3\hat{u}_{j\B})
   +\frac{2\B}{a^*}\inte\Big(Re(\rho_{1\B}^2\rho_{2\B}Q^3\hat{u}_{2\B})
                            +Re(\rho_{1\B}\rho_{2\B}^2Q^3\hat{u}_{1\B})\Big)\\
&=\sum^2_{j=1}Re\Big(\inte\Big\{
      \Big(\frac{\eps_\B^4\Omega^2}{4}|x|^2
           +\eps_\B^2\big(|\eps_\B x+x_\B|^2-1\big)^2\Big)\rho_{j\B}^2Q^2\\
&\qquad\qquad\qquad\quad+\Big(\frac{\eps_\B^4\Omega^2}{4}|x|^2
           +\eps_\B^2\big(|\eps_\B x+x_\B|^2-1\big)^2\Big)\hat{u}_{j\B}\rho_{j\B}Q\\
&\qquad\qquad\qquad\quad-\frac{a_j}{a^*}|\hat{u}_{j\B}|^2\hat{u}_{j\B}\rho_{j\B}Q
                        -\frac{2a_j}{a^*}\big(Re(\rho_{j\B}Q\hat{u}_{j\B})\big)^2
                        -\frac{a_j}{a^*}|\hat{u}_{j\B}|^2\rho_{j\B}^2Q^2
           \Big\}\Big)\\
&\quad-Re\Big(\inte\Big\{\frac{\B}{a^*}|\hat{u}_{2\B}|^2\hat{u}_{1\B}\rho_{1\B}Q
   +\frac{4\B}{a^*}Re\big(\rho_{1\B}\rho_{2\B}Q^2\hat{u}_{1\B}\hat{u}_{2\B}\big)
   +\frac{\B}{a^*}|\hat{u}_{2\B}|^2\rho_{1\B}^2Q^2\\
&\qquad\qquad\qquad\quad+\frac{\B}{a^*}|\hat{u}_{1\B}|^2\hat{u}_{2\B}\rho_{2\B}Q
   +\frac{\B}{a^*}|\hat{u}_{1\B}|^2\rho_{2\B}^2Q^2
           \Big\}\Big),
\end{split}
\end{equation}
where we have used the fact that
\begin{equation*}
  Re\Big(\inte-i\eps_\B^2\Omega \big(x^{\perp}\cdot\nabla\hat{u}_{j\B}\big)\rho_{j\B}Q\Big)
  =Re\Big(\inte i\eps_\B^2\Omega \big(x^{\perp}\cdot\nabla(\rho_{j\B}Q)\big)\hat{u}_{j\B}\Big)
  =0, \ \ j=1,2.
\end{equation*}
By Lemmas \ref{lem2.2} and \ref{lem2.3}, we deduce from \eqref{1.23:rhojB}, \eqref{2.53} and \eqref{2.55} that as $\B\nearrow\B^*$,
\begin{equation}\label{2.56}
\begin{split}
&\quad L.\ H.\ S.\ \ \hbox{of}\ \ \eqref{2.52}\\
&=\frac{(a^*)^2(2\B-a_1-a_2)}{\B^2-a_1a_2}-a^*\\
&\quad+\sum^2_{j=1}Re\Big(\inte\Big\{
      \Big(\frac{\eps_\B^4\Omega^2}{4}|x|^2
           +\eps_\B^2\big(|\eps_\B x+x_\B|^2-1\big)^2\Big)\rho_{j\B}^2Q^2\\
&\qquad\qquad\qquad\qquad+\Big(\frac{\eps_\B^4\Omega^2}{4}|x|^2
           +\eps_\B^2\big(|\eps_\B x+x_\B|^2-1\big)^2\Big)\rho_{j\B}Q\hat{u}_{j\B}
           \Big\}\Big)
+o(\eps_\B^{4})\\
&=\frac{(a^*)^2(2\B-a_1-a_2)}{\B^2-a_1a_2}-a^*\\
&\quad+\sum^2_{j=1}Re\Big(\inte
      \Big(\frac{\eps_\B^4\Omega^2}{4}|x|^2
      +\eps_\B^2\big(\eps_\B^2|x|^2+2\eps_\B x_2p_\B+p_\B^2-1\big)^2\Big)\rho_{j\B}^2Q^2\Big)\\
&\quad+\sum^2_{j=1}Re\Big(\inte
      \Big(\frac{\eps_\B^4\Omega^2}{4}|x|^2
      +\eps_\B^2\big(\eps_\B^2|x|^2+2\eps_\B x_2p_\B+p_\B^2-1\big)^2\Big)\rho_{j\B}^2Q\big(\eps_\B^4\psi_{1}+o(\eps_\B^4)\big)
      \Big)\\
&\quad+o(\eps_\B^{4})\\
&=\frac{(a^*)^2(2\B-a_1-a_2)}{\B^2-a_1a_2}-a^*
+\frac{a^*(2\B-a_1-a_2)}{\B^2-a_1a_2}\inte\frac{\eps_\B^4(\Omega^2+8)}{4}|x|^2Q^2
+o(\eps_\B^{4}).
\end{split}
\end{equation}

We now estimate the right hand side of \eqref{2.52}. By Lemma \ref{lem2.3}, we obtain that
\begin{equation}\label{2.57}
\begin{aligned}
  &\quad\sum^2_{j=1}\inte Re\l[i\eps_\B^2\Om(x^\perp\cdot\nabla v_{j\B})
                                         (x\cdot\nabla\ovl{v_{j\B}})\r]\\
  &=\sum^2_{j=1}\rho_{j\B}^2Re\Big(i\eps_\B^2\inte\l[x^\perp\cdot\nabla \l(Q+\eps_\B^4\psi_{1}+o(\eps_\B^4)\r)\r]
  \l[x\cdot\nabla \l(Q+\eps_\B^4\psi_{1}+o(\eps_\B^4)\r)\r]\Big)\\
  &=o(\eps_\B^{4})\ \ \hbox{as}\ \ \B\nearrow \B^*.
\end{aligned}
\end{equation}
Applying \eqref{1.23:rhojB} and Lemmas \ref{lem2.1}--\ref{lem2.3},
we derive from \eqref{2.16:exp} and \eqref{2.57} that
\begin{equation}\label{2.58}
\begin{split}
&\quad R.\ H.\ S.\ \ \hbox{of}\ \ \eqref{2.52}\\
&=\sum^2_{j=1}\inte\Big(\frac{\eps_\B^4\Omega^2}{2}|x|^2
                +\eps_\B^2\big(|\eps_\B x+x_\B|^2-1\big)^2
  +2\eps_\B^3(|\eps_\B x+x_\B|^2-1)
                (\eps_\B|x|^2+x\cdot x_\B)\Big)|v_{j\B}|^2\\
&\quad+o(\eps_\B^{4})\\
&=\sum^2_{j=1}Re\Big(
\inte\Big\{\Big[\frac{\eps_\B^4\Omega^2}{2}|x|^2
+\eps_\B^2\big(\eps_\B^2|x|^2+2\eps_\B x_2p_\B+p_\B^2-1\big)^2\Big]
\cdot\big[\rho_{j\B}^2Q^2+2\rho_{j\B}^2Q\eps_\B^4\psi_{1}+o(\eps_\B^4)\big]
      \Big\}\Big)\\
&\quad+2\eps_\B^3\sum^2_{j=1}Re\Big(\inte\Big\{
       \big[\eps_\B^3|x|^4+\l(3\eps_\B^2|x|^2+p_\B^2-1\r)x_{2}p_\B
            +2\eps_\B x_{2}^2 p_\B^2+\eps_\B(p_\B^2-1)|x|^2\big]\\
&\quad\qquad\qquad\qquad\qquad
\cdot\big[\rho_{j\B}^2Q^2+2\rho_{j\B}^2Q\eps_\B^4\psi_{1}+o(\eps_\B^4)\big]\Big\}\Big)
+o(\eps_\B^{4})\\
&=\frac{a^*(2\B-a_1-a_2)}{\B^2-a_1a_2}
\inte\frac{\eps_\B^4(\Omega^2+8)}{2}|x|^2Q^2+o(\eps_\B^{4})
\ \ \hbox{as}\ \ \B\nearrow\B^*.
\end{split}
\end{equation}

Combining \eqref{2.56} and \eqref{2.58} then yields that
\begin{equation*}
\begin{split}
&\quad\frac{(a^*)^2(2\B-a_1-a_2)}{\B^2-a_1a_2}-a^*\\
=&\frac{a^*(2\B-a_1-a_2)}{\B^2-a_1a_2}\eps_\B^4
\inte\frac{(\Omega^2+8)}{4}|x|^2Q^2
      +o(\eps_\B^{4})
\ \ \hbox{as}\ \ \B\nearrow\B^*,
\end{split}
\end{equation*}
from which we deduce that
\begin{equation}\label{2.59}
\frac{a^*-\frac{\B^2-a_1a_2}{2\B-a_1-a_2}}
{\eps_\B^4\frac{\Om^2+8}{4}\inte|x|^2Q^2}=1+o(1)
\ \ \hbox{as}\ \ \B\nearrow\B^*.
\end{equation}
Denote $F(\B):=\ds\frac{\B^2-a_1a_2}{2\B-a_1-a_2}$.
Using Taylor's expansion, we derive that
\begin{equation}\label{2.60}
\begin{split}
F(\B)&=F(\B^*)+F'(\B^*)(\B-\B^*)+o(|\B-\B^*|)\\
     &=a^*+2\gamma_1\gamma_2(\B-\B^*)
       +o(|\B-\B^*|) \ \ \hbox{as}\ \ \B\nearrow\B^*,\\
\end{split}
\end{equation}
where $\B^*>0$ is as in \eqref{1.11:B*}, and $\gam_j\in(0,1)$ for $j=1,2$ is defined by \eqref{1.15:gam}.
It then follows from \eqref{2.59} and \eqref{2.60} that
\begin{equation}\label{2.61}
\frac{2\gamma_1\gamma_2(\B^*-\B)+o(|\B-\B^*|)}
{\frac{(\Om^2+8)\eps_\B^4}{4}\inte|x|^2Q^2}
=1+o(1)
\ \ \hbox{as}\ \ \B\nearrow\B^*.
\end{equation}

Set $\alpha_\B
:=\Big(\frac{8\gamma_1\gamma_2(\B^*-\B)}{(\Om^2+8)\inte|x|^2Q^2}\Big)^\frac{1}{4}$, it then follows from \eqref{2.61} that \eqref{2.51:epal} holds true, and the proof of
Lemma \ref{lem2.4} is therefore complete.
\qed

\vskip 0.1truein

We are now already to complete the proof of Theorem \ref{thm1}.

\noindent {\bf Proof of Theorem \ref{thm1}.}
Using Lemma \ref{lem2.4}, we derive from \eqref{2.9:vjB} and \eqref{2.12:lim} that for $j=1,2$,
\begin{equation}\label{2.62:lim}
\sqrt{a^*}\alpha_\B u_{j\B}\big(\alpha_\B x+x_\B\big)
    e^{-i\big(\frac{\alpha_\B\Om}{2}x\cdot x_\B^\perp-\theta_{j\B}\big)}
\to\sqrt{\gamma_j}Q(x)\ \ \hbox{as}\ \ \B\nearrow\B^*
\end{equation}
strongly in $H^1(\R^2,\C)\cap L^\infty(\R^2,\C)$, where $\alpha_\B>0$ is defined by \eqref{1.16:alp}, $x_\B=(0,p_\B)$ is the unique maximal point of $|u_{1\B}|^2+|u_{2\B}|^2$, $\theta_{j\B}$ is chosen such that \eqref{2.10:theta} holds, and $\gam_j\in(0,1)$ is given by \eqref{1.15:gam}.
In view of \eqref{2.62:lim}, to prove \eqref{1.17:lim}, i.e., for $j=1,2$,
\begin{equation*}
\widetilde{v}_{j\B}(x):=\sqrt{a^*}\alpha_\B u_{j\B}(\alpha_\B x+x_\B)
e^{-i\big(\frac{\alpha_\B\Om}{2}x\cdot x_\B^\perp-\widetilde{\theta}_{j\B}\big)}
\to\sqrt{\gamma_j}Q(x)\ \ \hbox{as}\ \ \B\nearrow\B^*
\end{equation*}
strongly in $H^1(\R^2,\C)\cap L^\infty(\R^2,\C)$, we only need to prove that
$\lim\limits_{\B\nearrow\B^*}|\widetilde{\theta}_{j\B}-\theta_{j\B}|=0,\ \ j=1,2.$

Define
\begin{equation}\label{2.64:tvjB}
\widetilde{v}_{j\B}(x):=\sqrt{a^*}\alpha_\B u_{j\B}(\alpha_\B x+x_\B)
e^{-i\big(\frac{\alpha_\B\Om}{2}x\cdot x_\B^\perp-\widetilde{\theta}_{j\B}\big)},
\ \ j=1,2,
\end{equation}
where $\alpha_\B>0$ is defined by \eqref{1.16:alp}, $x_\B=(0,p_\B)$ is the unique maximal point of $|u_{1\B}|^2+|u_{2\B}|^2$, and $\widetilde{\theta}_{j\B}\in[0,2\pi)$ is chosen such that
\begin{equation}\label{2.65:theta}
\big\|\widetilde{v}_{j\B}-\sqrt{\gamma_j}Q\big\|_{L^2(\R^2)}
=\min\limits_{\theta\in[0,2\pi)}\big\|
\sqrt{a^*}\alpha_\B u_{j\B}(\alpha_\B x+x_\B)
e^{-i\big(\frac{\alpha_\B\Om}{2}x\cdot x_\B^\perp-\theta\big)}
-\sqrt{\gamma_j}Q\big\|_{L^2(\R^2)}.
\end{equation}
We then derive from \eqref{2.62:lim}, \eqref{2.64:tvjB} and \eqref{2.65:theta} that for $j=1,2$,
\begin{equation*}
\begin{aligned}
&\quad\lim_{\B\nearrow\B^*}
\big\|\sqrt{a^*}\alpha_\B u_{j\B}(\alpha_\B x+x_\B)
     e^{-i\big(\frac{\alpha_\B\Om}{2}x\cdot x_\B^\perp\big)}
            \big( e^{i\widetilde{\theta}_{j\B}}- e^{i\theta_{j\B}}\big)\big\|_{L^2(\R^2)}\\
&\le \lim_{\B\nearrow\B^*}
\big\|\sqrt{a^*}\alpha_\B u_{j\B}(\alpha_\B x+x_\B)
e^{-i\big(\frac{\alpha_\B\Om}{2}x\cdot x_\B^\perp-\widetilde{\theta}_{j\B}\big)}
-\sqrt{\gamma_j}Q\big\|_{L^2(\R^2)}\\
&\qquad+\lim_{\B\nearrow\B^*}
\big\|\sqrt{a^*}\alpha_\B u_{j\B}(\alpha_\B x+x_\B)
e^{-i\big(\frac{\alpha_\B\Om}{2}x\cdot x_\B^\perp-\theta_{j\B}\big)}
-\sqrt{\gamma_j}Q\big\|_{L^2(\R^2)}\\
&\le 2\lim_{\B\nearrow\B^*}
\big\|\sqrt{a^*}\alpha_\B u_{j\B}(\alpha_\B x+x_\B)
e^{-i\big(\frac{\alpha_\B\Om}{2}x\cdot x_\B^\perp-\theta_{j\B}\big)}
-\sqrt{\gamma_j}Q\big\|_{L^2(\R^2)}=0
\quad\hbox{as}\ \ \B\nearrow\B^*,
\end{aligned}
\end{equation*}
which yields that
\begin{equation}\label{2.66}
\lim_{\B\nearrow\B^*}|\widetilde{\theta}_{j\B}-\theta_{j\B}|=0,\ \ j=1,2.
\end{equation}
It then follows from \eqref{2.62:lim} and \eqref{2.66} that \eqref{1.17:lim} holds true. Furthermore, we obtain from \eqref{2.13:xB} that \eqref{1.18:xB} holds true.
The proof of Theorem \ref{thm1} is therefore complete.
\qed

\section{Axial Symmetry of Minimizers}
Suppose $V(x)$ satisfies \eqref{1.4:Vx}, and assume $0<\Om<\infty$,
$0< a_1,a_2<a^*$ are fixed. Following the limiting behavior of minimizers for $e(\Om,\f)$ established in Section 2, the main purpose of this section is to complete the proof of Theorem \ref{thm2} on the axial symmetry of minimizers for $e(\Om,\f)$.
We shall first establish the following proposition on the local uniqueness of minimizers for $e(\Om,\f)$ when $\B^*-\B>0$ is small enough.
\begin{prop}\label{prop3.1}
Suppose $V(x)$ satisfies \eqref{1.4:Vx}, and assume $0<\Om<\infty$ and
$0< a_1,a_2<a^*$ are fixed. If $(u_{1,1\B},u_{1,2\B})$ and $(u_{2,1\B},u_{2,2\B})$ are two complex-valued minimizers of $e(\Om,\f)$. Then we have
$$
\big(u_{1,1\B}(x),u_{1,2\B}(x)\big)\equiv
\big(u_{2,1\B}(\mathcal{R}x)e^{i\phi_{1\B}},u_{2,2\B}(\mathcal{R}x)e^{i\phi_{2\B}}\big)$$
for some constant phase $(\phi_{1\B},\phi_{2\B})\in[0,2\pi)\times[0,2\pi)$ and a suitable rotation $\mathcal{R}$ from $\R^2$ to $\R^2$ when $\B^*-\B>0$ is small enough.
\end{prop}

In order to prove Proposition \ref{prop3.1},
by contradiction, suppose that there exist two different minimizers $(u_{1,1\B},u_{1,2\B})$ and $(u_{2,1\B},u_{2,2\B})$ of $e(\Om,\f)$ as $\B\nearrow\B^*$.
Without loss of generality, we may assume
 \begin{equation}\label{3.1:assu}
   u_{1,1\B}(x)\not\equiv u_{2,1\B}(\mathcal{R}x)e^{i\psi_{1\B}}
 \ \ \hbox{in}\ \  \R^2
\end{equation}
for any rotation $\mathcal{R}:\R^2\to\R^2$ and the constant phase $\psi_{1\B}\in[0,2\pi)$.
Let $x_{j,\B}\in\R^2$ be the unique maximal point of $|u_{j,1\B}|^2+|u_{j,2\B}|^2$ for $j=1,2$. Note that the system \eqref{1.1:GPs} is invariant under the rotational transformation. Without loss of generality, we may assume $x_{j,\B}$ lies on the $x_2$-axis, i.e., $x_{j,\B}=(0,p_{j,\B})$.

Similar to the arguments of Lemmas \ref{lem2.2} and \ref{lem2.4}, one can deduce that
\begin{equation}\label{3.2}
p_{j,\B}-1=A\alpha_\B^2+o(\alpha_\B^2)
\ \ \hbox{as}\ \ \B\nearrow\B^*, \ \ j=1,2,
\end{equation}
where $A<0$ is as in \eqref{1.22:A}. It then follows that
\begin{equation}\label{3.3}
|x_{1,\B}-x_{2,\B}|=|p_{1,\B}-p_{2,\B}|=o(\alpha_\B^2)
\ \ \hbox{as}\ \ \B\nearrow\B^*.
\end{equation}
Define
\begin{equation}\label{3.4}
\left\{\begin{array}{lll}
\widetilde{v}_{j,1\B}(x):=\sqrt{a^*}\alpha_\B u_{j,1\B}(\alpha_\B x+x_{1,\B})
e^{-i\big(\frac{\alpha_\B\Om}{2}x\cdot
   x_{1,\B}^\perp-\widetilde{\theta}_{j,1\B}\big)},
\\[3mm]
\widetilde{v}_{j,2\B}(x):=\sqrt{a^*}\alpha_\B u_{j,2\B}(\alpha_\B x+x_{1,\B})
e^{-i\big(\frac{\alpha_\B\Om}{2}x\cdot
   x_{1,\B}^\perp-\widetilde{\theta}_{j,2\B}\big)},
\end{array}\right.
   \ \ j=1,2,
\end{equation}
where $\alpha_\B>0$ is as in \eqref{1.16:alp}, and
the constant phase
$(\widetilde{\theta}_{j,1\B},\widetilde{\theta}_{j,2\B})\in [0,2\pi)\times[0,2\pi)$ is chosen such that
\begin{equation}\label{3.5}
 Re\Big(\inte \widetilde{v}_{j,1\B}(iQ)\Big)
=Re\Big(\inte \widetilde{v}_{j,2\B}(iQ)\Big)=0, \  \  j=1,2.
\end{equation}
Using \eqref{3.3}, we obtain from Theorem \ref{thm1} that $\big(\widetilde{v}_{j,1\B}(x),\widetilde{v}_{j,2\B}(x)\big)$ satisfies for $\ j=1,2$,
\begin{equation}\label{3.6}
\left\{\begin{array}{lll}
\widetilde{v}_{j,1\B}(x)
\to\sqrt{\gamma_1}Q(x)\\[3mm]
\widetilde{v}_{j,2\B}(x)
\to\sqrt{\gamma_2}Q(x)
\end{array}\right.
\ \ \hbox{strongly in $H^1(\R^2,\C)\cap L^\infty(\R^2,\C)$ as}\ \ \B\nearrow\B^*,
\end{equation}
where $\gamma_1>0$ and $\gamma_2>0$ are defined by \eqref{1.15:gam}.

Denote
\begin{equation}\label{3.7}
\left\{\begin{array}{lll}
w_{1,1\B}(x):=\widetilde v_{1,1\B}(x)\!=\!R_{1,1\B}(x)+iI_{1,1\B}(x),\\[3mm]
w_{1,2\B}(x):=\widetilde v_{1,2\B}(x)\!=\!R_{1,2\B}(x)+iI_{1,2\B}(x),\\[3mm]
w_{2,1\B}(x):=\sqrt{a^*}\alpha_\B u_{2,1\B}
              \big(\mathcal{R}_0(\alpha_\B x+x_{1,\B})\big)
e^{-i\big(\frac{\alpha_\B\Om}{2}x\cdot x_{1,\B}^\perp-\theta_{2,1\B}\big)}
\!=\!R_{2,1\B}(x)+iI_{2,1\B}(x),
\\[3mm]
w_{2,2\B}(x):=\sqrt{a^*}\alpha_\B u_{2,2\B}
              \big(\mathcal{R}_0(\alpha_\B x+x_{1,\B})\big)
e^{-i\big(\frac{\alpha_\B\Om}{2}x\cdot x_{1,\B}^\perp-\theta_{2,2\B}\big)}
\!=\!R_{2,2\B}(x)+iI_{2,2\B}(x),
\end{array}\right.
\end{equation}
where $\mathcal{R}_0:\R^2\to\R^2$ is a suitable rotation such that $\mathcal{R}_0(x_{1,\B})$ satisfies
\begin{equation}\label{3.8}
  \big|u_{2,1\B}(\mathcal{R}_0(x_{1,\B}))\big|^2
  +\big|u_{2,2\B}(\mathcal{R}_0(x_{1,\B}))\big|^2
  =\max\limits_{|x|=|x_{1,\B}|}\big(|u_{2,1\B}(x)|^2+|u_{2,2\B}(x)|^2\big),
\end{equation}
and $(\theta_{2,1\B},\theta_{2,2\B})\in [0,2\pi)\times[0,2\pi)$ is chosen such that
\begin{equation}\label{3.9}
 Re\Big(\inte w_{2,1\B}(iQ)\Big)
=Re\Big(\inte w_{2,2\B}(iQ)\Big)=0.
\end{equation}
Here $\big(R_{j,1\B},R_{j,2\B}\big)$ and $\big(I_{j,1\B},I_{j,2\B}\big)$ denote the real and imaginary parts of $\big(w_{j,1\B},w_{j,2\B}\big)$ for $j=1,2$, respectively.
Following \eqref{3.4}, \eqref{3.7} and \eqref{3.8}, we deduce that
\begin{equation}\label{3.10}
  \frac{\partial\big(|w_{j,1\B}(x)|^2+|w_{j,2\B}(x)|^2\big)}
       {\partial x_{1}}\Big|_{x=0}=0,\ \ j=1,2.
\end{equation}
Similar to the argument of \cite [Lemma A.2] {GLLP}, one can derive from \eqref{3.3}--\eqref{3.9} that
\begin{equation}\label{3.11}
\left\{\begin{array}{lll}
w_{2,1\B}(x)
\to\sqrt{\gamma_1}Q(x)\\[3mm]
w_{2,2\B}(x)
\to\sqrt{\gamma_2}Q(x)
\end{array}\right.
\ \ \hbox{strongly in $H^1(\R^2,\C)\cap L^\infty(\R^2,\C)$ as}\ \ \B\nearrow\B^*.
\end{equation}
Combining \eqref{3.6}, \eqref{3.7} and \eqref{3.11} then yields that $\big(w_{j,1\B}(x),w_{j,2\B}(x)\big)$ satisfies for $\ j=1,2$,
\begin{equation}\label{3.12}
\left\{\begin{array}{lll}
w_{j,1\B}(x)
\to\sqrt{\gamma_1}Q(x)\\[3mm]
w_{j,2\B}(x)
\to\sqrt{\gamma_2}Q(x)
\end{array}\right.
\ \ \hbox{strongly in $H^1(\R^2,\C)\cap L^\infty(\R^2,\C)$ as}\ \ \B\nearrow\B^*.
\end{equation}
Note from \eqref{1.12:ELs}, \eqref{1.13:mu} and \eqref{3.7}
that for $j=1,2$, $\big(w_{j,1\B},w_{j,2\B}\big)$ satisfies the following system
\begin{equation}\label{3.13}
\left\{\begin{array}{lll}
  \ds\quad \mathcal{L}_\B w_{j,1\B}(x)
=\alpha^2_\B\mu_{j\B}w_{j,1\B}(x)
 +\frac{a_1}{a^*}|w_{j,1\B}|^2w_{j,1\B}(x)
 +\frac{\B}{a^*}|w_{j,2\B}|^2w_{j,1\B}(x)
\quad \hbox{in}\  \ \R^2,\\[5mm]
  \ds\quad \mathcal{L}_\B w_{j,2\B}(x)
=\alpha^2_\B\mu_{j\B}w_{j,2\B}(x)
 +\frac{a_2}{a^*}|w_{j,2\B}|^2w_{j,2\B}(x)
 +\frac{\B}{a^*}|w_{j,1\B}|^2w_{j,2\B}(x)
\quad\hbox{in}\  \ \R^2,
\end{array}\right.
\end{equation}
where the operator $\mathcal{L}_\B$ is defined by
\begin{equation}\label{3.14:LB}
\mathcal{L}_\B:
=-\Delta+i\,\alpha^2_\B\Om\,(x^\perp\cdot\nabla)
 +\frac{\alpha_\B^4\Omega^2}{4}|x|^2+\alpha_\B^2
 \big(|\alpha_\B x+x_{1,\B}|^2-1\big)^2
\ \ \hbox{in} \ \ \R^2,
\end{equation}
and $\mu_{j\B}\in\R$ satisfies
\begin{equation}\label{3.15:mu}
\mu_{j\B}=e(\Om,\f)-\frac{1}{(a^*)^2\alpha_\B^2}\inte \Big(\frac{a_1}{2}|w_{j,1\B}|^4
+\frac{a_2}{2}|w_{j,2\B}|^4
+\B|w_{j,1\B}|^2|w_{j,2\B}|^2\Big)dx.
\end{equation}
Similar to the argument of Lemma \ref{lem2.4}, one can deduce that the term $\mu_{j\B}\alpha_\B^2$ in \eqref{3.13} satisfies
\begin{equation}\label{3.16}
\alpha_\B^2\mu_{j\B}\to-1
\ \ \hbox{as}\ \ \B\nearrow\B^*, \ \ j=1,2.
\end{equation}
By the comparison principle,
one can derive from \eqref{3.12}, \eqref{3.13} and \eqref{3.16} that for $j,l=1,2$,
\begin{equation}\label{3.17}
      |w_{j,l\B}(x)|\leq Ce^{-\frac{2}{3}|x|} \quad \text{and} \quad
      |\nabla w_{j,l\B}(x)|\leq Ce^{-\frac{1}{2}|x|}
      \  \ \text{uniformly in} \  \  \R^2\  \ \text{as}\  \ \B\nearrow\B^*,
     \end{equation}
where the constant $C>0$ is independent of $0<\B<\B^*$.
Moreover, by the same argument of proving (5.10) in \cite{GGLL}, one can deduce from \eqref{3.12}--\eqref{3.17} that
\begin{equation}\label{3.18}
     C_1\|w_{2,2\B}-w_{1,2\B}\|_{L^\infty(\R^2)}
   \leq \|w_{2,1\B}-w_{1,1\B}\|_{L^\infty(\R^2)}
\leq C_2\|w_{2,2\B}-w_{1,2\B}\|_{L^\infty(\R^2)}
\  \ \text{as}\  \ \B\nearrow\B^*,
     \end{equation}
where the constants $C_1,C_2>0$ are independent of $\B$.
Using \eqref{2.2:dec}, \eqref{3.5}, \eqref{3.9}, \eqref{3.12}, \eqref{3.16} and \eqref{3.17}, the similar argument of \cite [Lemma 2.3] {GLP} yields from \eqref{3.13} that the imaginary parts $I_{j,l\B}$ of $w_{j,l\B}$ satisfies
for $j,l=1,2$,
\begin{equation}\label{3.19:IjBdecay}
      |I_{j,l\B}(x)|\leq C_{jl}(\alpha_\B)e^{-\frac{1}{4}|x|}, \quad
      |\nabla I_{j,l\B}(x)|\leq C_{jl}(\alpha_\B)e^{-\frac{1}{8}|x|} \ \
      \text{in} \  \  \R^2\ \ \hbox{as} \ \ \B\nearrow\B^*,
     \end{equation}
where the constant $C_{jl}(\alpha_\B)>0$ satisfies
$C_{jl}(\alpha_\B)=o(\alpha_\B^2)$ as $\B\nearrow\B^*$.

Under the assumption \eqref{3.1:assu}, we obtain from \eqref{3.4}, \eqref{3.7} and \eqref{3.18} that
\begin{equation}\label{3.20}
 w_{1,1\B}\not\equiv w_{2,1\B}\ \ \hbox{and}\ \
 w_{1,2\B}\not\equiv w_{2,2\B}\ \ \hbox{in}\ \ \R^2.
\end{equation}
Applying \eqref{3.12} and \eqref{3.16}--\eqref{3.20}, in the following we shall complete the proof of Proposition \ref{prop3.1}.

\vskip 0.1truein
\noindent {\bf Proof of Proposition \ref{prop3.1}.}
In view of \eqref{3.20}, we define
\begin{equation}\label{3.21}
 \left\{\begin{array}{lll}
  \xi_{1\B}(x):=\ds\frac{w_{2,1\B}(x)-w_{1,1\B}(x)}
  {\|w_{2,1\B}-w_{1,1\B}\|^{\frac{1}{2}}_{L^\infty(\R^2)}
   \|w_{2,2\B}-w_{1,2\B}\|^{\frac{1}{2}}_{L^\infty(\R^2)}}
   =R_{\xi_{1\B}}(x)+iI_{\xi_{1\B}}(x),\\[6mm]
 \xi_{2\B}(x):=\ds\frac{w_{2,2\B}(x)-w_{1,2\B}(x)}
  {\|w_{2,1\B}-w_{1,1\B}\|^{\frac{1}{2}}_{L^\infty(\R^2)}
   \|w_{2,2\B}-w_{1,2\B}\|^{\frac{1}{2}}_{L^\infty(\R^2)}}
   =R_{\xi_{2\B}}(x)+iI_{\xi_{2\B}}(x).
   \end{array}\right.
\end{equation}
We then derive from \eqref{3.13} that $\big(\xi_{1\B},\xi_{2\B}\big)$ satisfies
\begin{equation}\label{3.22}
 \left\{\begin{array}{lll}
\mathcal{L}_\B\xi_{1\B}
=\alpha^2_\B\mu_{2\B}\xi_{1\B}
+\ds\frac{\alpha^2_\B(\mu_{2\B}-\mu_{1\B})}
 {\|w_{2,1\B}-w_{1,1\B}\|^{\frac{1}{2}}_{L^\infty(\R^2)}
  \|w_{2,2\B}-w_{1,2\B}\|^{\frac{1}{2}}_{L^\infty(\R^2)}}
  w_{1,1\B}\\[5mm]
\qquad\qquad+\ds\frac{a_1}{a^*}|w_{2,1\B}|^2\xi_{1\B}
            +\frac{a_1}{a^*}\big(w_{2,1\B}\overline{\xi_{1\B}}
                                +\xi_{1\B}\overline{w_{1,1\B}}\big)w_{1,1\B}\\[3mm]
\qquad\qquad+\ds\frac{\B}{a^*}|w_{2,2\B}|^2\xi_{1\B}
            +\frac{\B}{a^*}\big(w_{2,2\B}\overline{\xi_{2\B}}
                                +\xi_{2\B}\overline{w_{1,2\B}}\big)w_{1,1\B}
\ \ \hbox{in} \ \ \R^2,\\[5mm]
\mathcal{L}_\B\xi_{2\B}
=\alpha^2_\B\mu_{2\B}\xi_{2\B}
+\ds\frac{\alpha^2_\B(\mu_{2\B}-\mu_{1\B})}
 {\|w_{2,1\B}-w_{1,1\B}\|^{\frac{1}{2}}_{L^\infty(\R^2)}
  \|w_{2,2\B}-w_{1,2\B}\|^{\frac{1}{2}}_{L^\infty(\R^2)}}
  w_{1,2\B}\\[5mm]
\qquad\qquad+\ds\frac{a_1}{a^*}|w_{2,2\B}|^2\xi_{2\B}
            +\frac{a_1}{a^*}\big(w_{2,2\B}\overline{\xi_{2\B}}
                                +\xi_{2\B}\overline{w_{1,2\B}}\big)w_{1,2\B}\\[3mm]
\qquad\qquad+\ds\frac{\B}{a^*}|w_{2,1\B}|^2\xi_{2\B}
            +\frac{\B}{a^*}\big(w_{2,1\B}\overline{\xi_{1\B}}
                                +\xi_{1\B}\overline{w_{1,1\B}}\big)w_{1,2\B}
\ \ \hbox{in} \ \ \R^2,
\end{array}\right.
\end{equation}
where the operator $\mathcal{L}_\B$ is given by \eqref{3.14:LB}.
Furthermore, we obtain from \eqref{3.18} that $\big(\xi_{1\B},\xi_{2\B}\big)$ satisfies
\begin{equation}\label{3.23}
0\leq|\xi_{1\B}(x)|,\ |\xi_{2\B}(x)|\leq C<\infty\quad \text{and}
\quad |\xi_{1\B}(x)\xi_{2\B}(x)|\leq1\  \ \text{in}\  \ \R^2,
\end{equation}
where the constant $C>0$ is independent of $\B$.

We are now ready to carry out the proof of Proposition \ref{prop3.1} by the following four steps:
\vskip 0.1truein

{\em\noindent $\mathbf{Step\  \ 1}$. Passing to a subsequence if necessary, the function $\big(\xi_{1\B},\xi_{2\B}\big)$ defined by \eqref{3.21} satisfies
\begin{equation}\label{3.24}
\xi_{j\B}(x)\to
\sqrt{\gamma_j}
\Big(c_0(Q+x\cdot\nabla Q)+b_0^j(iQ)+\sum^2_{l=1}b_l\frac{\partial Q}{\partial x_l}\Big)
\ \ \hbox{in}\ \ C^{1}_{loc}(\R^2,\C)\ \ \hbox{as}\ \ \B\nearrow\B^*,
\end{equation}
where $j=1,2$, $\gam_j\in(0,1)$ is as in \eqref{1.15:gam}, $c_0$, $b_0^j$, $b_1$ and $b_2$ are some constants independent of $\B$.}

We first claim that
\begin{equation}\label{3.25}
\frac{|\alpha^2_\B(\mu_{2\B}-\mu_{1\B})|}
 {\|w_{2,1\B}-w_{1,1\B}\|^{\frac{1}{2}}_{L^\infty(\R^2)}
  \|w_{2,2\B}-w_{1,2\B}\|^{\frac{1}{2}}_{L^\infty(\R^2)}}\leq C
\ \ \hbox{as}\ \ \B\nearrow\B^*,
\end{equation}
where the constant $C>0$ is independent of $\B$.
Indeed, we obtain from \eqref{3.15:mu} that
\begin{equation}\label{3.26}
\begin{split}
&\quad\frac{\alpha^2_\B(\mu_{2\B}-\mu_{1\B})}
  {\|w_{2,1\B}-w_{1,1\B}\|^{\frac{1}{2}}_{L^\infty(\R^2)}
  \|w_{2,2\B}-w_{1,2\B}\|^{\frac{1}{2}}_{L^\infty(\R^2)}}\\
  &=-Re\Big(\inte\Big[\frac{a_1}{2(a^*)^2}
\big(\big|w_{2,1\B}\big|^2+\big|w_{1,1\B}\big|^2\big)
\big(w_{2,1\B}+w_{1,1\B}\big)\overline{\xi_{1\B}}\\
&\qquad\qquad\qquad+\frac{a_2}{2(a^*)^2}
\big(\big|w_{2,2\B}\big|^2+\big|w_{1,2\B}\big|^2\big)
\big(w_{2,2\B}+w_{1,2\B}\big)\overline{\xi_{2\B}}
\Big]\Big)\\
&\quad-\frac{\B}{(a^*)^2}Re\Big(
\inte\Big[\big|w_{2,2\B}\big|^2
\big(w_{2,1\B}+w_{1,1\B}\big)\overline{\xi_{1\B}}
+\big|w_{1,1\B}\big|^2
\big(w_{2,2\B}+w_{1,2\B}\big)\overline{\xi_{2\B}}
\Big]\Big).
\end{split}
\end{equation}
Using \eqref{3.17} and \eqref{3.23}, we deduce from \eqref{3.26} that \eqref{3.25} holds true.

By the same argument of proving (5.15) in \cite{GGLL}, one can deduce from  \eqref{3.16}, \eqref{3.17}, \eqref{3.22}, \eqref{3.23} and \eqref{3.25} that for $j=1,2$,
\begin{equation}\label{3.27}
|\xi_{j\B}(x)|\leq Ce^{-\frac{2}{3}|x|} \quad \text{and} \quad
|\nabla \xi_{j\B}(x)|\leq Ce^{-\frac{1}{2}|x|}
\ \ \text{uniformly in} \  \  \R^2\  \ \text{as}\  \ \B\nearrow\B^*,
\end{equation}
where the constant $C>0$ is independent of $0<\B<\B^*$.
It follows from \eqref{3.27} that $(x^\perp\cdot\nabla\xi_{j\B})$ is bounded uniformly and decays exponentially for sufficiently large $|x|$ as $\B\nearrow\B^*$.
Applying \eqref{3.16}, \eqref{3.17}, \eqref{3.25} and \eqref{3.27},
the standard elliptic regularity theory then yields from \eqref{3.22} that $\|\xi_{j\B}\|_{C^{1,\alpha}_{loc}(\R^2,\C)}\leq C$ uniformly as $\B\nearrow\B^*$ for some $\alpha\in(0,1)$. Therefore, up to a subsequence if necessary,
 \begin{equation}\label{3.28}
 \xi_{j\B}\to\xi_{j0}\ \ \hbox{in}\ \ C^{1}_{loc}(\R^2,\C)
 \ \ \hbox{as} \ \ \B\nearrow\B^*, \ \ j=1,2.
 \end{equation}
Using \eqref{3.12}, \eqref{3.27} and \eqref{3.28}, we derive from \eqref{3.26} that
\begin{equation}\label{3.29}
\begin{split}
&\quad\lim_{\B\nearrow\B^*}
 \frac{\alpha^2_\B(\mu_{2\B}-\mu_{1\B})}
  {\|w_{2,1\B}-w_{1,1\B}\|^{\frac{1}{2}}_{L^\infty(\R^2)}
  \|w_{2,2\B}-w_{1,2\B}\|^{\frac{1}{2}}_{L^\infty(\R^2)}}\\
&=-\frac{1}{a^*}Re\Big\{\Big(\frac{2a_1\gamma_1^{\frac{3}{2}}}{a^*}
           +\frac{2\B^*\gamma_2\sqrt{\gamma_1}}{a^*}\Big)
           \inte Q^3\overline{\xi_{10}}
+\Big(\frac{2a_2\gamma_2^{\frac{3}{2}}}{a^*}
           +\frac{2\B^*\gamma_1\sqrt{\gamma_2}}{a^*}\Big)
           \inte Q^3\overline{\xi_{20}}
\Big\}\\
&=-Re\Big(\frac{2\sqrt{\gamma_1}}{a^*}\inte Q^3\overline{\xi_{10}}
         +\frac{2\sqrt{\gamma_2}}{a^*}\inte Q^3\overline{\xi_{20}}
\Big),
\end{split}
\end{equation}
where $\gam_1,\gam_2\in(0,1)$ are as in \eqref{1.15:gam}, and we have used the facts that $a_1\gamma_1+\B^*\gamma_2=a^*$ and $a_2\gamma_2+\B^*\gamma_1=a^*$.
By \eqref{3.12}, \eqref{3.16} and \eqref{3.27}--\eqref{3.29}, we deduce from \eqref{3.22} that  $(\xi_{10},\xi_{20})$ solves the following system
\begin{equation}\label{3.30}
 \left\{\begin{array}{lll}

\mathcal{L}_{1}(\xi_{10},\xi_{20})
=\ds-\sqrt{\gamma_1}QRe\Big(\frac{2\sqrt{\gamma_1}}{a^*}\inte Q^3\overline{\xi_{10}}
         +\frac{2\sqrt{\gamma_2}}{a^*}\inte Q^3\overline{\xi_{20}}\Big)
\quad \hbox{in}\  \ \R^2,\\[5mm]
\mathcal{L}_{2}(\xi_{20},\xi_{10}) =\ds-\sqrt{\gamma_2}QRe\Big(\frac{2\sqrt{\gamma_2}}{a^*}\inte Q^3\overline{\xi_{20}}
         +\frac{2\sqrt{\gamma_1}}{a^*}\inte Q^3\overline{\xi_{10}}\Big)
\quad \hbox{in}\  \ \R^2,
\end{array}\right.
\end{equation}
where $\mathcal{L}_{1}$ and $\mathcal{L}_{2}$ are defined by \eqref{2.29}.
Following \cite [Lemma 2.2 and Theorem 3.1] {DW}, we obtain that the solution of the following linearized system
\begin{equation}\label{3.31:linop}
 \left\{\begin{array}{lll}
 \begin{split}
\mathcal{L}_{1}(\phi_1,\phi_2)=0
\  \ \hbox{in}\  \ \R^2,\\[3mm]
\mathcal{L}_{2}(\phi_2,\phi_1)=0
\  \ \hbox{in}\  \ \R^2,
 \end{split}
\end{array}\right.
\end{equation}
is given by
\begin{equation}\label{3.32:nond}
\begin{pmatrix}
  \phi_1\\[5mm]
  \phi_2
\end{pmatrix}
=\begin{pmatrix}
ib_0^1\sqrt{\gamma_1}Q\\[5mm]
ib_0^2\sqrt{\gamma_2}Q
\end{pmatrix}
+
\sum^2_{l=1}b_l
\begin{pmatrix}
\ds\sqrt{\gamma_1}\frac{\partial Q}{\partial x_l}\\[5mm]
\ds\sqrt{\gamma_2}\frac{\partial Q}{\partial x_l}
\end{pmatrix}
\end{equation}
for some constants $b_0^j,b_l\in\R$, where $j,l=1,2$.
Moreover, one can check that
\begin{equation}\label{3.33}
 \left\{\begin{array}{lll}
\mathcal{L}_1
\big(\sqrt{\gamma_1}(Q+x\cdot\nabla Q),\sqrt{\gamma_2}(Q+x\cdot\nabla Q)\big)
=-2\sqrt{\gamma_1}Q,\\[3mm]
\mathcal{L}_2
\big(\sqrt{\gamma_2}(Q+x\cdot\nabla Q),\sqrt{\gamma_1}(Q+x\cdot\nabla Q)\big)
=-2\sqrt{\gamma_2}Q.
\end{array}\right.
\end{equation}
We then derive from \eqref{3.30}--\eqref{3.33} that there exist constants $c_0$, $b_0^j$ and $b_l$, where $j,l=1,2$, such that
\begin{equation*}
\xi_{j0}=
\sqrt{\gamma_j}
\Big(c_0(Q+x\cdot\nabla Q)+b_0^j(iQ)+\sum^2_{l=1}b_l\frac{\partial Q}{\partial x_l}\Big)
\ \ \hbox{in}\ \ \R^2, \ \ j=1,2,
\end{equation*}
together with \eqref{3.28} then yields that \eqref{3.24} holds true.

\vskip 0.1truein
{\em\noindent $\mathbf{Step\  \ 2}$. The constants $b_0^1$, $b_0^2$ and $b_2$ in \eqref{3.24} satisfy $b_0^1=b_0^2=b_2=0$.}

Note from \eqref{3.5}, \eqref{3.7} and \eqref{3.9} that $Re\big(\inte\xi_{j\B}(iQ)\big)=0$ for $j=1,2$. We then derive from \eqref{3.24} and \eqref{3.27} that
\begin{equation*}
  Re\Big(\inte\sqrt{\gamma_j}
\Big(c_0(Q+x\cdot\nabla Q)+b_0^j(iQ)+\sum^2_{l=1}b_l\frac{\partial Q}{\partial x_l}\Big)(iQ)\Big)=0, \ \ j=1,2,
\end{equation*}
which implies that $b_0^j=0$ for $j=1,2$.

We next prove that $b_2=0$. Multiplying the first equation of \eqref{3.13} by $\frac{\p\ovl{w_{j,1\B}}}{\p x_2}
+\frac{-i\alpha_\B^2\Om x_{1}\ovl{w_{j,1\B}}}{2}$ and multiplying the second equation of \eqref{3.13} by
$\frac{\p\ovl{w_{j,2\B}}}{\p x_2}
+\frac{-i\alpha_\B^2\Om x_{1}\ovl{w_{j,2\B}}}{2}$,
integrating over $\R^2$ and taking its real part, the same argument of proving
\cite [Lemma A.1] {GGLL} shows that
\begin{equation}\label{3.34}
\inte\frac{\p \big(|\alpha_\B x+x_{1,\B}|^2-1\big)^2}{\p x_{2}}
\big(|w_{j,1\B}|^2+|w_{j,2\B}|^2\big)=0
\ \ \hbox{as}\ \ \B\nearrow\B^*, \ \ j=1,2.
\end{equation}
By \eqref{3.2}, \eqref{3.12} and \eqref{3.24}, we derive from \eqref{3.34} that
\begin{equation*}
\begin{aligned}
0&=Re\Big(\inte\frac{\p \big(|\alpha_\B x+x_{1,\B}|^2-1\big)^2}{\p x_2}
       \Big[(w_{2,1\B}+w_{1,1\B})\overline{\xi_{1\B}}
           +(w_{2,2\B}+w_{1,2\B})\overline{\xi_{2\B}}\Big]\Big)\\
&=Re\Big(\inte 4\alpha_\B\big(\big|\alpha_\B x+x_{1,\B}\big|^2-1\big)
\big(\alpha_\B x_2+p_{1,\B}\big)
\Big[(w_{2,1\B}+w_{1,1\B})\overline{\xi_{1\B}}
    +(w_{2,2\B}+w_{1,2\B})\overline{\xi_{2\B}}\Big]\Big)\\
 &=Re\Big(\inte 4\alpha_\B
  \big(\alpha_\B^2|x|^2+p_{1,\B}^2-1+2\alpha_\B x_2 p_{1,\B}\big)
  \big(\alpha_\B x_2+p_{1,\B}\big)\\
&\qquad\qquad\qquad\quad\Big[(w_{2,1\B}+w_{1,1\B})\overline{\xi_{1\B}}
                 +(w_{2,2\B}+w_{1,2\B})\overline{\xi_{2\B}}\Big]\Big)\\
 &=Re\Big(\inte 8\alpha_\B^2 x_2
 \Big[(w_{2,1\B}+w_{1,1\B})\overline{\xi_{1\B}}
    +(w_{2,2\B}+w_{1,2\B})\overline{\xi_{2\B}}\Big]
 \Big)+o(\alpha_\B^2)\\
&=16\alpha_\B^2(\gamma_1+\gamma_2)
  \inte x_2Q\Big[c_{0}(Q+x\cdot\nabla Q)+\sum_{l=1}^2b_l\frac{\p Q}{\p x_l}\Big]
  +o(\alpha_\B^2)\\
&=-8\alpha_\B^2b_2\inte Q^2+o(\alpha_\B^2)
\ \ \hbox{as}\ \ \B\nearrow\B^*,
\end{aligned}
\end{equation*}
which implies that $b_2=0$.

\vskip 0.1truein
{\em\noindent $\mathbf{Step\  \ 3}$. The constant $c_0$ in \eqref{3.24} satisfies $c_0=0$.}

Similar to the argument of \eqref{2.52}, one can deduce from \eqref{3.13} that $\big(w_{j,1\B},w_{j,2\B}\big)$ satisfies the following Pohozaev identity
\begin{equation}\label{3.35}
\begin{split}
&\quad Re\Big(\inte
\Big[i\alpha_\B^2\Om(x^\perp\cdot\nabla w_{j,1\B})(x\cdot\nabla\ovl{w_{j,1\B}})
+i\alpha_\B^2\Om(x^\perp\cdot\nabla w_{j,2\B})(x\cdot\nabla\ovl{w_{j,2\B}})\Big]
\Big)\\
&=-\alpha_\B^2\mu_{j\B}\inte\big(|w_{j,1\B}|^2+|w_{j,2\B}|^2\big)
+\alpha_\B^2\inte\Big(
\frac{\alpha_\B^2\Omega^2}{2}|x|^2+\big(|\alpha_\B x+x_{1,\B}|^2-1\big)^2\\
&\qquad+2\alpha_\B(|\alpha_\B x+x_{1,\B}|^2-1)(\alpha_\B|x|^2+x\cdot x_{1,\B})
\Big)\big(|w_{j,1\B}|^2+|w_{j,2\B}|^2\big)\\
&\quad-\frac{a_1}{2a^*}\inte|w_{j,1\B}|^4-\frac{a_2}{2a^*}\inte|w_{j,2\B}|^4
-\frac{\B}{a^*}\inte|w_{j,1\B}|^2|w_{j,2\B}|^2,\ \ j=1,2.
\end{split}
\end{equation}
We then derive from \eqref{3.35} that
\begin{equation}\label{3.36}
\begin{split}
&\quad\frac{A_{2\B}-A_{1\B}}
{\|w_{2,1\B}-w_{1,1\B}\|^{\frac{1}{2}}_{L^\infty(\R^2)}
  \|w_{2,2\B}-w_{1,2\B}\|^{\frac{1}{2}}_{L^\infty(\R^2)}}\\
  &=\frac{(B_{2\B}-B_{1\B})+(C_{2\B}-C_{1\B})}
{\|w_{2,1\B}-w_{1,1\B}\|^{\frac{1}{2}}_{L^\infty(\R^2)}
  \|w_{2,2\B}-w_{1,2\B}\|^{\frac{1}{2}}_{L^\infty(\R^2)}},
\end{split}
\end{equation}
where $A_{j\B}$, $B_{j\B}$ and $C_{j\B}$ for $j=1,2$ are defined by
\begin{equation*}
A_{j\B}:=Re\Big(\inte \Big[i\,\alpha^2_\B\Om\big(x^\perp\cdot\nabla w_{j,1\B}\big)
                        \big(x\cdot\nabla \overline{w_{j,1\B}}\big)
                          +i\,\alpha^2_\B\Om\big(x^\perp\cdot\nabla w_{j,2\B}\big)
                        \big(x\cdot\nabla \overline{w_{j,2\B}}\big)
\Big]\Big),
\end{equation*}
\begin{equation*}
\begin{split}
B_{j\B}:=&-\alpha_\B^2\mu_{j\B}\inte\big(|w_{j,1\B}|^2+|w_{j,2\B}|^2\big)
+\alpha_\B^2\inte\Big(
\frac{\alpha_\B^2\Omega^2}{2}|x|^2+\big(|\alpha_\B x+x_{1,\B}|^2-1\big)^2\\
&\qquad+2\alpha_\B(|\alpha_\B x+x_{1,\B}|^2-1)(\alpha_\B|x|^2+x\cdot x_{1,\B})
\Big)\big(|w_{j,1\B}|^2+|w_{j,2\B}|^2\big),
\end{split}
\end{equation*}
and
\begin{equation*}
C_{j\B}:=-\frac{a_1}{2a^*}\inte|w_{j,1\B}|^4
         -\frac{a_2}{2a^*}\inte|w_{j,2\B}|^4
         -\frac{\B}{a^*}\inte|w_{j,1\B}|^2|w_{j,2\B}|^2.
\end{equation*}

We are now going to estimate all the terms in \eqref{3.36}. By \eqref{3.16}, \eqref{3.17}, \eqref{3.19:IjBdecay} and \eqref{3.25}, the similar argument of \cite [Proposition 3.4] {GLP} yields from \eqref{3.22} that
\begin{equation}\label{3.37}
\|\nabla I_{\xi_{j\B}}(x)\|_{L^\infty(\R^2)}=O(\alpha_\B^2)
\ \ \hbox{as} \ \ \B\nearrow\B^*, \ \ j=1,2.
\end{equation}
Using \eqref{3.7} and \eqref{3.21}, we deduce that
\begin{equation}\label{3.38}
\begin{split}
&\quad\frac{A_{2\B}-A_{1\B}}
{\|w_{2,1\B}-w_{1,1\B}\|^{\frac{1}{2}}_{L^\infty(\R^2)}
  \|w_{2,2\B}-w_{1,2\B}\|^{\frac{1}{2}}_{L^\infty(\R^2)}}\\
&=Re\Big(\inte \Big[i\,\alpha^2_\B\Om\big(x^\perp\cdot\nabla \xi_{1\B}\big)
                        \big(x\cdot\nabla \overline{w_{2,1\B}}\big)
                          +i\,\alpha^2_\B\Om\big(x^\perp\cdot\nabla w_{1,1\B}\big)
                        \big(x\cdot\nabla \overline{\xi_{1\B}}\big)\\
&\qquad\qquad\quad+i\,\alpha^2_\B\Om\big(x^\perp\cdot\nabla \xi_{2\B}\big)
                        \big(x\cdot\nabla \overline{w_{2,2\B}}\big)
                          +i\,\alpha^2_\B\Om\big(x^\perp\cdot\nabla w_{1,2\B}\big)
                        \big(x\cdot\nabla \overline{\xi_{2\B}}\big)
\Big]\Big)\\
&=\inte \alpha^2_\B\Om\big(x^\perp\cdot\nabla R_{\xi_{1\B}}\big)
                       \big(x\cdot\nabla I_{2,1\B}\big)
-\inte \alpha^2_\B\Om\big(x^\perp\cdot\nabla I_{\xi_{1\B}}\big)
                       \big(x\cdot\nabla R_{2,1\B}\big)\\
&\quad+\inte \alpha^2_\B\Om\big(x^\perp\cdot\nabla R_{1,1\B}\big)
                       \big(x\cdot\nabla I_{\xi_{1\B}}\big)
-\inte \alpha^2_\B\Om\big(x^\perp\cdot\nabla I_{1,1\B}\big)
                       \big(x\cdot\nabla R_{\xi_{1\B}}\big)\\
&\quad+\inte \alpha^2_\B\Om\big(x^\perp\cdot\nabla R_{\xi_{2\B}}\big)
                       \big(x\cdot\nabla I_{2,2\B}\big)
-\inte \alpha^2_\B\Om\big(x^\perp\cdot\nabla I_{\xi_{2\B}}\big)
                       \big(x\cdot\nabla R_{2,2\B}\big)\\
&\quad+\inte \alpha^2_\B\Om\big(x^\perp\cdot\nabla R_{1,2\B}\big)
                       \big(x\cdot\nabla I_{\xi_{2\B}}\big)
-\inte \alpha^2_\B\Om\big(x^\perp\cdot\nabla I_{1,2\B}\big)
                       \big(x\cdot\nabla R_{\xi_{2\B}}\big)\\
&=
\inte \alpha^2_\B\Om\big[x^\perp\cdot\nabla
                   \big(x\cdot\nabla(\sqrt{\gamma_1}Q)\big) \big]I_{\xi_{1\B}}
+\inte \alpha^2_\B\Om\big[x^\perp\cdot\nabla
                   \big(x\cdot\nabla(\sqrt{\gamma_2}Q)\big) \big]I_{\xi_{2\B}}
+o(\alpha^4_\B)\\
&=o(\alpha^4_\B)\ \ \hbox{as}\ \ \B\nearrow\B^*,
\end{split}
\end{equation}
where we have used \eqref{3.12}, \eqref{3.19:IjBdecay}, \eqref{3.27} and \eqref{3.37} in the third equality, and the fact that $Q(x)$ is radially symmetric in the last equality.
As for the term containing $B_{j\B}$, applying \eqref{3.2}, we derive from \eqref{3.12}, \eqref{3.24}, \eqref{3.27} and Step 2 that
\begin{equation}\label{3.39}
\begin{split}
&\quad\frac{B_{2\B}-B_{1\B}}
 {\|w_{2,1\B}-w_{1,1\B}\|^{\frac{1}{2}}_{L^\infty(\R^2)}
  \|w_{2,2\B}-w_{1,2\B}\|^{\frac{1}{2}}_{L^\infty(\R^2)}}\\
&=-\frac{\alpha^2_\B(\mu_{2\B}-\mu_{1\B})a^*}
       {\|w_{2,1\B}-w_{1,1\B}\|^{\frac{1}{2}}_{L^\infty(\R^2)}
        \|w_{2,2\B}-w_{1,2\B}\|^{\frac{1}{2}}_{L^\infty(\R^2)}}\\
&\quad+\alpha_\B^2Re\Big(\inte\Big\{
\frac{\alpha_\B^2\Omega^2}{2}|x|^2+\big(|\alpha_\B x+x_{1,\B}|^2-1\big)^2\\
&\qquad\qquad\qquad\qquad
+2\alpha_\B(|\alpha_\B x+x_{1,\B}|^2-1)(\alpha_\B|x|^2+x\cdot x_{1,\B})
\Big\}\\
&\qquad\qquad\qquad
\cdot\Big[
\big(\xi_{1\B}\overline{w_{2,1\B}}+w_{1,1\B}\overline{\xi_{1\B}}\big)
+\big(\xi_{2\B}\overline{w_{2,2\B}}+w_{1,2\B}\overline{\xi_{2\B}}\big)
\Big]\Big)\\
&=-\frac{\alpha^2_\B(\mu_{2\B}-\mu_{1\B})a^*}
       {\|w_{2,1\B}-w_{1,1\B}\|^{\frac{1}{2}}_{L^\infty(\R^2)}
        \|w_{2,2\B}-w_{1,2\B}\|^{\frac{1}{2}}_{L^\infty(\R^2)}}\\
&\quad+2(\gamma_1+\gamma_2)\alpha_\B^2\inte
\Big(\frac{\alpha_\B^2\Omega^2}{2}|x|^2+8\alpha_\B^2x_2^2\Big)
\Big[Q\Big(c_0(Q+x\cdot\nabla Q)+b_1\frac{\p Q}{\p x_1}\Big)\Big]
+o(\alpha_\B^4)\\
&=-\frac{\alpha^2_\B(\mu_{2\B}-\mu_{1\B})a^*}
       {\|w_{2,1\B}-w_{1,1\B}\|^{\frac{1}{2}}_{L^\infty(\R^2)}
        \|w_{2,2\B}-w_{1,2\B}\|^{\frac{1}{2}}_{L^\infty(\R^2)}}\\
&\quad-\alpha_\B^4(\Omega^2+8)c_0\inte |x|^2Q^2+o(\alpha_\B^4)
\ \ \hbox{as}\ \ \B\nearrow\B^*,
\end{split}
\end{equation}
where $\gamma_j$ is as in \eqref{1.15:gam} for $j=1,2$, and we have used the facts that $\rho_{j\B}\to\sqrt{\gam_j}$ as $\B\nearrow\B^*$, and $Q(x)$ is radially symmetric. For the term containing $C_{j\B}$, we deduce from \eqref{3.26} that
\begin{equation}\label{3.40}
\begin{split}
&\quad\frac{C_{2\B}-C_{1\B}}
  {\|w_{2,1\B}-w_{1,1\B}\|^{\frac{1}{2}}_{L^\infty(\R^2)}
  \|w_{2,2\B}-w_{1,2\B}\|^{\frac{1}{2}}_{L^\infty(\R^2)}}\\
  &=-Re\Big(\inte\Big[\frac{a_1}{2a^*}
\big(\big|w_{2,1\B}\big|^2+\big|w_{1,1\B}\big|^2\big)
\big(w_{2,1\B}+w_{1,1\B}\big)\overline{\xi_{1\B}}\\
&\qquad\qquad\qquad+\frac{a_2}{2a^*}
\big(\big|w_{2,2\B}\big|^2+\big|w_{1,2\B}\big|^2\big)
\big(w_{2,2\B}+w_{1,2\B}\big)\overline{\xi_{2\B}}
\Big]\Big)\\
&\quad-\frac{\B}{a^*}Re\Big(
\inte\Big[\big|w_{2,2\B}\big|^2
\big(w_{2,1\B}+w_{1,1\B}\big)\overline{\xi_{1\B}}
+\big|w_{1,1\B}\big|^2
\big(w_{2,2\B}+w_{1,2\B}\big)\overline{\xi_{2\B}}
\Big]\Big)\\
&=\frac{\alpha^2_\B(\mu_{2\B}-\mu_{1\B})a^*}
  {\|w_{2,1\B}-w_{1,1\B}\|^{\frac{1}{2}}_{L^\infty(\R^2)}
  \|w_{2,2\B}-w_{1,2\B}\|^{\frac{1}{2}}_{L^\infty(\R^2)}}.
\end{split}
\end{equation}

Following \eqref{3.38}--\eqref{3.40}, we obtain from \eqref{3.36} that
\begin{equation*}
\alpha_\B^4(\Omega^2+8)c_0\inte |x|^2Q^2=o(\alpha_\B^4)
\ \ \hbox{as}\ \ \B\nearrow\B^*,
\end{equation*}
which implies that $c_0=0$.

\vskip 0.1truein
{\em\noindent $\mathbf{Step\  \ 4}$. The constant $b_1$ in \eqref{3.24} satisfies $b_1=0$.}

Multiplying $\xi_{1\B}$ by $\ovl{w_{2,1\B}+w_{1,1\B}}$ and multiplying $\xi_{2\B}$ by $\ovl{w_{2,2\B}+w_{1,2\B}}$, one can deduce from \eqref{3.21} that
\begin{equation}\label{3.41}
\begin{aligned}
  &\quad\xi_{1\B}\big(\ovl{w_{2,1\B}+w_{1,1\B}}\big)
  +\xi_{2\B}\big(\ovl{w_{2,2\B}+w_{1,2\B}}\big)\\
  &=\frac{|w_{2,1\B}|^2+|w_{2,2\B}|^2-\big(|w_{1,1\B}|^2+|w_{1,2\B}|^2\big)}
  {\|w_{2,1\B}-w_{1,1\B}\|^{\frac{1}{2}}_{L^\infty(\R^2)}
  \|w_{2,2\B}-w_{1,2\B}\|^{\frac{1}{2}}_{L^\infty(\R^2)}}\\[3mm]
  &\quad+\xi_{1\B}\ovl{w_{1,1\B}}-w_{1,1\B}\ovl{\xi_{1\B}}
  +\xi_{2\B}\ovl{w_{1,2\B}}-w_{1,2\B}\ovl{\xi_{2\B}}
\ \ \hbox{in}\ \ \R^2.
\end{aligned}
\end{equation}
Denote $\partial_{1}f:=\frac{\partial f}{\partial x_1}$,
we then derive from \eqref{3.10} and \eqref{3.41} that
\begin{equation}\label{3.42}
\begin{aligned}
&\quad\partial_{1}\big[\xi_{1\B}\big(\ovl{w_{2,1\B}+w_{1,1\B}}\big)
                      +\xi_{2\B}\big(\ovl{w_{2,2\B}+w_{1,2\B}}\big)\big]\Big|_{x=0}\\
&=\partial_{1}\big[\xi_{1\B}\ovl{w_{1,1\B}}-w_{1,1\B}\ovl{\xi_{1\B}}
                    +\xi_{2\B}\ovl{w_{1,2\B}}-w_{1,2\B}\ovl{\xi_{2\B}}\big]\Big|_{x=0}\\
&=\Big[(\partial_{1}\xi_{1\B})\ovl{w_{1,1\B}}+\xi_{1\B}(\partial_{1}\ovl{w_{1,1\B}})
-(\partial_{1}\ovl{\xi_{1\B}})w_{1,1\B}-(\partial_{1}w_{1,1\B})\ovl{\xi_{1\B}}\\
&\qquad+(\partial_{1}\xi_{2\B})\ovl{w_{1,2\B}}+\xi_{2\B}(\partial_{1}\ovl{w_{1,2\B}})
-(\partial_{1}\ovl{\xi_{2\B}})w_{1,2\B}-(\partial_{1}w_{1,2\B})\ovl{\xi_{2\B}}\Big]
\Big|_{x=0}.
\end{aligned}
\end{equation}
By Steps 2 and 3, we obtain that the constants $c_0,b_0^1,b_0^2,b_2$ in \eqref{3.24} satisfies $c_0=b_0^1=b_0^2=b_2=0$.
We then deduce from \eqref{3.12} and \eqref{3.24} that
\begin{equation}\label{3.43}
  R.\ H.\ S.\ \ \hbox{of}\ \  \eqref{3.42}\ \ \to 0\ \ \hbox{as}\ \ \B\nearrow\B^*.
\end{equation}
On the other hand, using \eqref{3.12} and \eqref{3.24} again, we deduce that
\begin{equation}\label{3.44}
\begin{aligned}
  L.\ H.\ S.\ \ \hbox{of}\ \  \eqref{3.42}&=\partial_{1}\big[\xi_{1\B}\big(\ovl{w_{2,1\B}+w_{1,1\B}}\big)
                      +\xi_{2\B}\big(\ovl{w_{2,2\B}+w_{1,2\B}}\big)\big]\Big|_{x=0}\\
  &=\Big[\partial_{1}\xi_{1\B}\big(\ovl{w_{2,1\B}+w_{1,1\B}}\big)
         +\xi_{1\B}\partial_{1}\big(\ovl{w_{2,1\B}+w_{1,1\B}}\big)\\
         &\qquad+\partial_{1}\xi_{2\B}\big(\ovl{w_{2,2\B}+w_{1,2\B}}\big)
         +\xi_{2\B}\partial_{1}\big(\ovl{w_{2,2\B}+w_{1,2\B}}\big)
         \Big]\Big|_{x=0}\\
  &\to 2Q(0)b_1\frac{\p^2 Q(0)}{\p x_1^2}\ \ \hbox{as}\ \ \B\nearrow\B^*,
\end{aligned}
\end{equation}
where we have used the fact that $\rho_{j\B}\to\sqrt{\gam_j}$ as $\B\nearrow\B^*$ for $j=1,2$.
It then follows from \eqref{3.43} and \eqref{3.44} that $b_1=0$.

We now take $(x_\B,y_\B)$ such that $|\xi_{1\B}(x_\B)\xi_{2\B}(y_\B)|=\|\xi_{1\B}\xi_{2\B}\|_{L^\infty(\R^2)}=1$.
Applying the exponential decay of $\xi_{1\B}$ and $\xi_{2\B}$ in \eqref{3.27}, we deduce that $|x_\B|\leq C$ and $|y_\B|\leq C$ uniformly in $\B$. Furthermore, one can conclude that $|\xi_{j\B}|\to\xi_{j0}\not\equiv0$ uniformly in $C^1_{loc}(\R^2)$ as $\B\nearrow\B^*$, where $j=1,2$.
However, we obtain from Steps 1--4 that $\xi_{j0}\equiv0$, this is a contradiction, and hence the assumption \eqref{3.1:assu} is false. We then have
$$\big(u_{1,1\B}(x),u_{1,2\B}(x)\big)\equiv
\big(u_{2,1\B}(\mathcal{R}x)e^{i\phi_{1\B}},u_{2,2\B}(\mathcal{R}x)e^{i\phi_{2\B}}\big)$$
for some constant phase $(\phi_{1\B},\phi_{2\B})\in[0,2\pi)\times[0,2\pi)$ and a suitable rotation $\mathcal{R}$ from $\R^2$ to $\R^2$ when $\B^*-\B>0$ is small enough.
The proof of Proposition \ref{prop3.1} is therefore complete.
\qed

\vskip 0.1truein

Following Proposition \ref{prop3.1}, we now give the proof of Theorem \ref{thm2}.

\noindent {\bf Proof of Theorem \ref{thm2}.}
Let $(u_1,u_2)\in\m$, where $\m$ is as in \eqref{1.8:m}. Note that
\begin{equation}\label{3.45}
\begin{split}
\inte x^{\perp}\cdot\big(iu_j,\nabla u_j\big)
&=\inte Re\big(x^{\perp}\cdot iu_j\nabla \bar{u}_j\big)\\
&=-\inte Re\big(x^{\perp}\cdot i\bar{u}_j\nabla u_j\big)
=-\inte x^{\perp}\cdot\big(i\bar{u}_j,\nabla \bar{u}_j\big), \ \ j=1,2.
\end{split}
\end{equation}
Here and below $\bar{f}$ denotes the conjugate of $f$. We then derive from \eqref{1.7:GPf} and \eqref{3.45} that
\begin{equation*}
\begin{split}
 F_{-\Om,\f}(\bar{u}_1,\bar{u}_2)
 &=\sum_{j=1}^2\inte\Big[|\nabla \bar{u}_j|^2+V(x)|\bar{u}_j|^2
            -\frac{a_j}{2}|\bar{u}_j|^4
            +\Om\, x^{\perp}\cdot\big(i\bar{u}_j,\nabla \bar{u}_j\big)\Big]dx\\
 &\quad-\inte\B|\bar{u}_1|^2|\bar{u}_2|^2dx\\
 &=F_{\Om,\f}(u_1,u_2),
\end{split}
\end{equation*}
which further yields that
\begin{equation}\label{3.46}
 e(-\Om,\f)=e(\Om,\f).
\end{equation}

On the other hand, we note that
\begin{equation}\label{3.47}
\begin{split}
&\quad\inte x^{\perp}\cdot\big(iu_j(-x_1,x_2),\nabla u_j(-x_1,x_2)\big)\\
&=\inte Re\big(x^{\perp}\cdot iu_j(-x_1,x_2)\nabla \overline{u_j(-x_1,x_2)}\big)\\
&=-\inte Re\big(x^{\perp}\cdot i u_j(x_1,x_2)\nabla \overline{u_j(x_1,x_2)}\big)\\
&=-\inte x^{\perp}\cdot\big(iu_j(x_1,x_2),\nabla u_j(x_1,x_2)\big), \ \ j=1,2.
\end{split}
\end{equation}
By \eqref{3.47}, we derive from \eqref{1.7:GPf} that
\begin{equation}\label{3.48}
\begin{split}
 &\quad F_{-\Om,\f}\big(u_1(-x_1,x_2),u_2(-x_1,x_2)\big)\\
 &=\sum_{j=1}^2\inte\Big[|\nabla u_j(-x_1,x_2)|^2+V(x)|u_j(-x_1,x_2)|^2
            -\frac{a_j}{2}|u_j(-x_1,x_2)|^4\\
            &\qquad+\Om\, x^{\perp}\cdot\big(iu_j(-x_1,x_2),\nabla u_j(-x_1,x_2)\big)\Big]dx
 -\inte\B|u_1(-x_1,x_2)|^2|u_2(-x_1,x_2)|^2dx\\
 &=F_{\Om,\f}\big(u_1(x_1,x_2),u_2(x_1,x_2)\big).
\end{split}
\end{equation}

Let $(u_{1\B},u_{2\B})$ be a minimizer of $e(\Om,\f)$ and $x_\B=(0,p_\B)\in\R^2$ be the unique maximal point of $|u_{1\B}|^2+|u_{2\B}|^2$. Using \eqref{3.46} and \eqref{3.48}, we obtain that $(u_{1\B}(-x_1,x_2),u_{2\B}(-x_1,x_2))$
is also a minimizer of $e(\Om,\f)$. It then follows from Proposition \ref{prop3.1} that
\begin{equation}\label{3.49:axial}
\big(u_{1\B}(-x_1,x_2),u_{2\B}(-x_1,x_2)\big)\equiv
\big(u_{1\B}(\mathcal{R}x)e^{i\phi_{1\B}},
     u_{2\B}(\mathcal{R}x)e^{i\phi_{2\B}}\big)
\end{equation}
for some constant phase $(\phi_{1\B},\phi_{2\B})\in[0,2\pi)\times[0,2\pi)$ and a suitable rotation $\mathcal{R}$ from $\R^2$ to $\R^2$ when $\B^*-\B>0$ is small enough. Furthermore, one can derive from \eqref{3.49:axial} that $x_\B=(0,p_\B)$ is also the maximal point of $|u_{1\B}(-x_1,x_2)|^2+|u_{2\B}(-x_1,x_2)|^2$.
By the uniqueness of the maximal point $x_\B=(0,p_\B)$, we derive from \eqref{3.49:axial} that $\mathcal{R}$ is an identity transformation, i.e. $\mathcal{R}=id$. Therefore, we further obtain that
\begin{equation}\label{3.50}
\big(u_{1\B}(-x_1,x_2),u_{2\B}(-x_1,x_2)\big)\equiv
\big(u_{1\B}(x)e^{i\phi_{1\B}},
     u_{2\B}(x)e^{i\phi_{2\B}}\big)
\end{equation}
for some constant phase $(\phi_{1\B},\phi_{2\B})\in[0,2\pi)\times[0,2\pi)$ when $\B^*-\B>0$ is small enough.
In view of \eqref{3.50}, we obtain that Theorem \ref{thm2} holds true.
\qed

\section{Refined Spike Profiles of Minimizers}
Suppose $V(x)$ satisfies \eqref{1.4:Vx}, and assume $0<\Om<\infty$,
$0< a_1,a_2<a^*$ are fixed. Let $(u_{1\B},u_{2\B})$ be the unique minimizer of $e(\Om,\f)$  and $x_\B$ be the unique maximal point of $|u_{1\B}|^2+|u_{2\B}|^2$. Note that the system \eqref{1.1:GPs} is invariant under the rotational transformation. Without loss of generality, we may assume that $x_{\B}$ lies on the $x_2$-axis, i.e., $x_{\B}=(0,p_{\B})$.
The main purpose of this section is to establish Theorem \ref{thm3} on the refined spike profiles of minimizers $(u_{1\B},u_{2\B})$ for $e(\Om,\f)$ as $\B\nearrow\B^*$.

Let $v_{j\B}(x):=\sqrt{a^*}\eps_\B u_{j\B}\big(\eps_\B x+x_\B\big)
    e^{-i\big(\frac{\eps_\B\Om}{2}x\cdot x_\B^\perp-\theta_{j\B}\big)}$
be defined by \eqref{2.9:vjB}, where $\eps_\B$ is as in \eqref{2.7:eps}, $x_\B=(0,p_\B)\in\R^2$ is the unique maximal point of $|u_{1\B}|^2+|u_{2\B}|^2$, and the constant phase $\theta_{j\B}\in [0,2\pi)$ is chosen such that \eqref{2.10:theta} holds true.
We then obtain from  \eqref{2.16:exp}, Lemma \ref{lem2.3} and \eqref{2.45} that
\begin{equation}\label{4.1:expan}
v_{j\B}(x):=\rho_{j\B}Q(x)+\hat{u}_{j\B}(x)
=\rho_{j\B}Q(x)+\rho_{j\B}\eps_\B^4\psi_{1}(x)+o(\eps_\B^4)
\ \ \hbox{as} \ \ \B\nearrow\B^*,\ \ j=1,2,
\end{equation}
where $\rho_{j\B}>0$ is as in \eqref{1.23:rhojB}, $\psi_{1}(x)$ is the unique solution of \eqref{2.43:psi1}, and $\hat{u}_{j\B}(x)$ satisfies
\begin{equation}\label{4.2:lowdec}
|\hat{u}_{j\B}(x)|\leq C\eps_\B^4e^{-\frac{1}{8}|x|}, \ \
|\nabla\hat{u}_{j\B}(x)|\leq C\eps_\B^4e^{-\frac{1}{16}|x|}
\ \ \hbox{as}\ \ \B\nearrow\B^*, \ \ j=1,2,
\end{equation}
where the constant $C>0$ is independent of $0<\B<\B^*$.
Furthermore, we obtain from Lemmas \ref{lem2.2} and \ref{lem2.4} that the maximal point $x_\B=(0,p_\B)$ of $|u_{1\B}|^2+|u_{2\B}|^2$ satisfies
\begin{equation}\label{4.3:pB}
p_{\B}-1=A\eps_\B^2+o(\eps_\B^2)
\ \ \hbox{as}\ \ \B\nearrow\B^*,
\end{equation}
where $A=-\frac{\inte |x|^2Q^2dx}{a^*}<0$ is given by \eqref{1.22:A}
and
\begin{equation}\label{4.5:epal}
  \eps_\B=\alpha_\B+o(\alpha_\B)\ \ \hbox{as}\ \ \B\nearrow\B^*,
\end{equation}
where $\eps_\B$ and $\alpha_\B>0$ are as in \eqref{2.7:eps} and \eqref{1.16:alp}, respectively.

\subsection{Refined spike profiles of  $(v_{1\B},v_{2\B})$}
In this subsection, we shall derive the refined spike profiles of  $(v_{1\B},v_{2\B})$ defined in \eqref{2.9:vjB}.
We first give the more refined expansion of the maximal point $x_\B$ for $|u_{1\B}|^2+|u_{2\B}|^2$ than \eqref{4.3:pB} as $\B\nearrow\B^*$.
\begin{lem}\label{lem4.1}
Suppose $V(x)$ satisfies \eqref{1.4:Vx}, and assume $0<\Om<\infty$ and
$0< a_1,a_2<a^*$ are fixed.
Let $(u_{1\B},u_{2\B})$ be the unique minimizer of $e(\Om,\f)$ and $x_\B=(0,p_\B)$ be the unique maximal point of $|u_{1\B}|^2+|u_{2\B}|^2$. Then we have
\begin{equation}\label{4.6:pB}
p_{\B}-1=A\eps_\B^2-\frac{A^2}{2}\eps_\B^4+o(\eps_\B^4)
\ \ \hbox{as}\ \ \B\nearrow\B^*,
\end{equation}
where $A<0$ is as in \eqref{1.22:A}.
\end{lem}
\noindent{\bf Proof.}
Recall from \eqref{2.40} that when $\B$ close to $\B^*$ enough,
\begin{equation}\label{4.7}
\inte\frac{\p \big(|\eps_\B x+x_\B|^2-1\big)^2}
{\p x_{2}}\big(|v_{1\B}|^2+|v_{2\B}|^2\big)=0.
\end{equation}
Denote
\begin{equation}\label{4.8}
B_\B:=p_\B-1-A\eps_\B^2.
\end{equation}
By \eqref{4.1:expan}--\eqref{4.3:pB}, we derive from \eqref{4.7} and \eqref{4.8} that
\begin{equation*}
\begin{aligned}
0&=\sum^2_{j=1}\inte\frac{\p \big(|\eps_\B x+x_\B|^2-1\big)^2}
{\p x_{2}}|v_{j\B}|^2dx\\
&=\sum^2_{j=1}\inte4\eps_\B(|\eps_\B x+x_\B|^2-1)(\eps_\B x_2+p_\B)
                         \big|\rho_{j\B}Q+\hat{u}_{j\B}\big|^2dx\\
&=\sum^2_{j=1}\inte4\eps_\B(\eps_\B^2|x|^2+2\eps_\B x_2p_\B+p_\B^2-1)
                           (\eps_\B x_2+p_\B)
             \l(\rho_{j\B}^2Q^2+2Re(\rho_{j\B}Q\hat{u}_{j\B})+|\hat{u}_{j\B}|^2\r)dx\\
&=\sum^2_{j=1}\rho_{j\B}^2\Big[
  4\eps_\B^3\inte(|x|^2+2A+2x_{2}^2)Q^2dx
 +4\eps_\B^5\inte(3A^2+A|x|^2+2Ax_{2}^2)Q^2dx\\
 &\qquad\qquad\quad+8\eps_\B B_\B\inte Q^2dx\Big]
  +o(\eps_\B^5)\\
&=\sum^2_{j=1}\rho_{j\B}^2\big(4\eps_\B^5A^2a^*+8\eps_\B B_\B a^*\big)+o(\eps_\B^5)
\ \ \hbox{as}\ \ \B\nearrow\B^*,
\end{aligned}
\end{equation*}
which yields that
\begin{equation}\label{4.9}
  B_\B=-\frac{A^2}{2}\eps_\B^4+o(\eps_\B^4)\ \ \hbox{as}\ \ \B\nearrow\B^*.
\end{equation}
It then follows from \eqref{4.8} and \eqref{4.9} that \eqref{4.6:pB} holds true, and
Lemma \ref{lem4.1} is thus proved.
\qed

\vskip 0.1truein

Based on Lemma \ref{lem4.1}, we are devoted to deriving the refined expansion of $\big(v_{1\B}(x),v_{2\B}(x)\big)$ in terms of $\eps_\B$.
\begin{lem}\label{lem4.2}
Under the assumptions of Lemma \ref{lem4.1}, let $v_{j\B}(x)$ be defined by \eqref{2.9:vjB}, where $j=1,2$. Then we have for $j=1,2$,
\begin{equation}\label{4.10}
\begin{aligned}
v_{j\B}(x):=&\rho_{j\B}Q(x)+\hat{u}_{j\B}(x)\\
=&\rho_{j\B}Q(x)+\rho_{j\B}\eps_\B^4\psi_{1}(x)+\rho_{j\B}\eps_\B^5\psi_{2}(x)
+\rho_{j\B}\eps_\B^6\psi_{3}(x)+o(\eps_\B^6)\ \ \hbox{in}\ \ \R^2
\end{aligned}
\end{equation}
as $\B\nearrow\B^*$, where $\rho_{j\B}>0$ is as in \eqref{1.23:rhojB}, $\psi_{k}(x)\in C^2(\R^2,\C)\cap L^\infty(\R^2,\C)$ solves uniquely
\begin{equation}\label{4.11:psik}
\mathcal{L}\psi_{k}(x)=f_k(x)\ \ \mbox{in}\ \ \R^2,
\ \ \nabla Re(\psi_{k}(0))=0,\ Re\inte \psi_k(iQ)=0,\ \ k=1,2,3,
\end{equation}
where $f_k(x)$ satisfies
\begin{equation*}\arraycolsep=1.5pt
f_k(x)=\left\{\begin{array}{lll}
&-\big(\frac{\Om^2}{4}|x|^2+4x_{2}^2\big)Q(x),\qquad \quad \ &\mbox{if}\ \ k=1;\\[3mm]
&-\big(4|x|^2+8A\big)x_{2}Q(x), \,\ &\mbox{if}\ \ k=2;\\[3mm]
&-\big(|x|^4+4A|x|^2+8Ax_{2}^2+4A^2\big)Q(x)-i\Om(x^\perp\cdot\nabla\psi_1)
, \,\ &\mbox{if}\ \ k=3.
\end{array}\right.
\end{equation*}
Here $A<0$ and the operator $\mathcal{L}$ are as in \eqref{1.22:A} and \eqref{2.3:L}, respectively.
\end{lem}
\noindent{\bf Proof.}
Set
\begin{equation}\label{4.13}
\widetilde{\mathcal{U}}_{j\B}(x)
:=\hat{u}_{j\B}(x)-\rho_{j\B}\eps_\B^4\psi_{1}(x)
   -\rho_{j\B}\eps_\B^5\psi_{2}(x)-\rho_{j\B}\eps_\B^{6}\psi_{3}(x),
\ \ j=1,2,
\end{equation}
where $\rho_{j\B}>0$ is defined by \eqref{1.23:rhojB}, and
$\psi_{k}(x)\in C^2(\R^2,\C)\cap L^\infty(\R^2,\C)$ is a solution of
\eqref{4.11:psik} for $k=1,2,3$. It then follows from \eqref{2.4:Lnon} and \eqref{4.11:psik} that $\psi_{k}(x)$ is unique for $k=1,2,3$.
Note from \eqref{2.18:GPs}, \eqref{2.20:LF} and \eqref{4.11:psik}--\eqref{4.13} that
\begin{equation}\label{4.14}
\left\{\begin{array}{lll}
\begin{split}
\mathcal{L}_{1\B}(\widetilde{\mathcal{U}}_{1\B},\widetilde{\mathcal{U}}_{2\B})
 &=F_{1\B}(x)+\rho_{1\B}\eps_\B^4\big(\frac{\Om^2}{4}|x|^2+4x_{2}^2\big)Q
 +\rho_{1\B}\eps_\B^5\big(4|x|^2+8A\big)x_{2}Q\\[3mm]
 &\quad+\rho_{1\B}\eps_\B^6
 \Big[\big(|x|^4+4A|x|^2+8Ax_{2}^2+4A^2\big)Q+i\,\Om(x^\perp\cdot\nabla\psi_1)\Big]\\[3mm]
&:=\widetilde{\mathcal{F}}_{1\B}(x)
\quad \mbox{in}\,\  \R^2,\\[3mm]
\mathcal{L}_{2\B}(\widetilde{\mathcal{U}}_{2\B},\widetilde{\mathcal{U}}_{1\B})
 &=F_{2\B}(x)+\rho_{2\B}\eps_\B^4\big(\frac{\Om^2}{4}|x|^2+4x_{2}^2\big)Q
 +\rho_{2\B}\eps_\B^5\big(4|x|^2+8A\big)x_{2}Q\\[3mm]
 &\quad+\rho_{2\B}\eps_\B^6
 \Big[\big(|x|^4+4A|x|^2+8Ax_{2}^2+4A^2\big)Q+i\,\Om(x^\perp\cdot\nabla\psi_1)\Big]\\[3mm]
&:=\widetilde{\mathcal{F}}_{2\B}(x)
\quad \mbox{in}\,\  \R^2,
\end{split}
\end{array}\right.
\end{equation}
where the operators $\mathcal{L}_{1\B}$ and $\mathcal{L}_{2\B}$ are as in \eqref{2.19:linop}, and $F_{j\B}(x)$ is given by \eqref{2.20:LF} for $j=1,2$.
Using \eqref{2.2:dec} and Lemma \ref{lem4.1}, we deduce from \eqref{4.1:expan} and \eqref{4.2:lowdec} that the term $\widetilde{\mathcal{F}}_{j\B}(x)$ satisfies for $j=1,2$,
\begin{equation}\label{4.15}
\begin{split}
 \ds|\widetilde{\mathcal{F}}_{j\B}(x)|:=&\Big|F_{j\B}(x)
                  +\rho_{j\B}\eps_\B^4\Big(\frac{\Om^2}{4}|x|^2+4x_{2}^2\Big)Q
                  +\rho_{j\B}\eps_\B^5\big(4|x|^2+8A\big)x_{2}Q\\[3mm]
                  \quad\qquad\qquad&+\rho_{j\B}\eps_\B^6
 \Big[\big(|x|^4+4A|x|^2+8Ax_{2}^2+4A^2\big)Q+i\,\Om(x^\perp\cdot\nabla\psi_1)\Big]
                  \Big|\\[3mm]
 \le& C_{j\B}\eps_\B^6e^{-\frac{1}{10}|x|}
 \ \ \hbox{uniformly in}\ \ \R^2 \ \ \hbox{as} \ \ \B\nearrow \B^*,\\[3mm]
\end{split}
\end{equation}
where $C_{j\B}>0$ satisfies $C_{j\B}=o(1)$ as $\B\nearrow\B^*$ for $j=1,2$.
Similar to the argument of \eqref{2.22}, one can deduce from \eqref{4.2:lowdec}, \eqref{4.14} and \eqref{4.15} that
\begin{equation}\label{4.16}
\big|\widetilde{\mathcal{U}}_{j\B}(x)\big|
\leq C_{j\B} \eps_\B^6\quad \hbox{uniformly in}\ \ \R^2
\ \ \hbox{as}\, \ \B\nearrow \B^*,\ \  j=1,2.
\end{equation}
We then conclude from \eqref{4.13} and \eqref{4.16} that \eqref{4.10} holds true, and the proof of Lemma \ref{lem4.2} is therefore complete.
\qed

\subsection{Proof of Theorem \ref{thm3}}
In this subsection, we shall complete the proof of Theorem \ref{thm3}. We shall first use the refined expansions of $x_\B=(0,p_\B)$ and $\big(v_{1\B}(x),v_{2\B}(x)\big)$ established in Lemmas \ref{lem4.1} and \ref{lem4.2} to derive the following refined expansions of $\eps_\B$ and $\mu_{\B}$ in terms of $\alpha_\B$.
\begin{lem}\label{lem4.3}
Under the assumptions of Lemma \ref{lem4.1}, the term $\eps_\B$ defined by \eqref{2.7:eps} satisfies
\begin{equation}\label{4.17}
  \eps_\B=\alpha_\B-\frac{B_2}{4B_1}\alpha_\B^3
  +\frac{9B_2^2-8B_1^2B_4}{32B_1^2}\alpha_\B^5+o(\alpha_\B^5)
  \ \ \hbox{as}\ \ \B\nearrow\B^*,
\end{equation}
where $\alpha_\B>0$ is given by \eqref{1.16:alp},
and $B_1,B_2,B_3,B_4$ are as in \eqref{1.21:B1234}.
\end{lem}
\noindent{\bf Proof.}
We shall use the Pohozaev identity \eqref{2.52} to prove \eqref{4.17}. For the left hand side of \eqref{2.52}, using \eqref{1.23:rhojB}, \eqref{4.2:lowdec}, Lemmas \ref{lem4.1} and \ref{lem4.2}, the similar argument of proving \eqref{2.56} yields that
\begin{equation}\label{4.18}
\begin{split}
&\quad L.\ H.\ S.\ \ \hbox{of}\ \ \eqref{2.52}\\
&=\frac{(a^*)^2(2\B-a_1-a_2)}{\B^2-a_1a_2}-a^*\\
&\quad+\frac{a^*(2\B-a_1-a_2)}{\B^2-a_1a_2}
\inte\Big[\frac{\eps_\B^4(\Omega^2+8)}{4}|x|^2Q^2
+\eps_\B^6\big(|x|^4+8A|x|^2+4A^2\big)Q^2\\
&\qquad\qquad\qquad\qquad\qquad\qquad
+\eps_\B^8\Big(\frac{\Om^2}{4}|x|^2+4x_2^2\Big)Q\psi_1
\Big]
+o(\eps_\B^{8})
\ \ \hbox{as}\ \ \B\nearrow\B^*.
\end{split}
\end{equation}
As for the right hand side of \eqref{2.52}, we first obtain from Lemma \ref{lem4.2} that
\begin{equation}\label{4.19}
\begin{aligned}
  &\quad\sum^2_{j=1}\inte Re\l[i\eps_\B^2\Om(x^\perp\cdot\nabla v_{j\B})
                                         (x\cdot\nabla\ovl{v_{j\B}})\r]\\
  &=\sum^2_{j=1}\rho_{j\B}^2Re\Big(i\eps_\B^2\inte\l[x^\perp\cdot\nabla \l(Q+\eps_\B^4\psi_{1}+\eps_\B^5\psi_{2}+\eps_\B^6\psi_{3}
  +o(\eps_\B^6)\r)\r]\\
  &\qquad\qquad\qquad\qquad\quad
  \cdot\l[x\cdot\nabla \l(Q+\eps_\B^4\psi_{1}+\eps_\B^5\psi_{2}+\eps_\B^6\psi_{3}
  +o(\eps_\B^6)\r)\r]\Big)\\
  &=o(\eps_\B^{8})\ \ \hbox{as}\ \ \B\nearrow \B^*,
\end{aligned}
\end{equation}
where $\psi_k$ are given by \eqref{4.11:psik} for $k=1,2,3$, and we have used the fact that
\begin{equation*}
  \inte(x^\perp\cdot\nabla\psi_{k})(x\cdot\nabla Q)
  =-\inte\psi_{k}x^\perp\cdot\nabla(x\cdot\nabla Q)=0.
\end{equation*}
By \eqref{1.23:rhojB}, \eqref{4.2:lowdec}, Lemma \ref{lem4.1} and \ref{lem4.2}, and \eqref{4.19}, the similar argument of proving \eqref{2.58} shows that
\begin{equation}\label{4.20}
\begin{split}
&\quad R.\ H.\ S.\ \ \hbox{of}\ \ \eqref{2.52}\\
&=\frac{a^*(2\B-a_1-a_2)}{\B^2-a_1a_2}
\inte\Big[\frac{\eps_\B^4(\Omega^2+8)}{2}|x|^2Q^2
+\eps_\B^6\big(3|x|^4+16A|x|^2+4A^2\big)Q^2\\
&\qquad\qquad\qquad\qquad\qquad\qquad
+\eps_\B^8\Big(\Om^2|x|^2+16x_2^2\Big)Q\psi_1
\Big]
+o(\eps_\B^{8})
\ \ \hbox{as}\ \ \B\nearrow\B^*.
\end{split}
\end{equation}

Combining \eqref{4.18} and \eqref{4.20} then yields that
\begin{equation*}
\begin{split}
\quad&\frac{(a^*)^2(2\B-a_1-a_2)}{\B^2-a_1a_2}-a^*\\
=&\frac{a^*(2\B-a_1-a_2)}{\B^2-a_1a_2}\inte\Big[\frac{\eps_\B^4(\Omega^2+8)}{4}|x|^2Q^2
+2\eps_\B^6\big(|x|^4+4A|x|^2\big)Q^2
+3\eps_\B^8\Big(\frac{\Om^2}{4}|x|^2+4x_2^2\Big)Q\psi_1
\Big]\\
&+o(\eps_\B^{8})\\
:=&\frac{a^*(2\B-a_1-a_2)}{\B^2-a_1a_2}
\big(B_1\eps_\B^4+B_2\eps_\B^6+B_3\eps_\B^8\big)+o(\eps_\B^{8})
\ \ \hbox{as}\ \ \B\nearrow\B^*,
\end{split}
\end{equation*}
from which we deduce that
\begin{equation}\label{4.21}
\frac{a^*-\frac{\B^2-a_1a_2}{2\B-a_1-a_2}}{B_1\eps_\B^4}
=1+\frac{B_2}{B_1}\eps_\B^2+\frac{B_3}{B_1}\eps_\B^4+o(\eps_\B^{4})
\ \ \hbox{as}\ \ \B\nearrow\B^*.
\end{equation}
Denote $F(\B):=\ds\frac{\B^2-a_1a_2}{2\B-a_1-a_2}$.
Using Taylor's expansion, we deduce that
\begin{equation*}
\begin{split}
F(\B)&=F(\B^*)+F'(\B^*)(\B-\B^*)+\frac{F''(\B^*)}{2}(\B-\B^*)^2+o(|\B-\B^*|^2)\\
     &=a^*-2\gamma_1\gamma_2(\B^*-\B)
       +\frac{1-4\gamma_1\gamma_2}{2\B^*-a_1-a_2}(\B^*-\B)^2+o(|\B^*-\B|^2)
        \ \ \hbox{as}\ \ \B\nearrow\B^*,
\end{split}
\end{equation*}
where $\B^*>0$ is as in \eqref{1.11:B*}, and $\gam_j\in(0,1)$ for $j=1,2$ is defined by \eqref{1.15:gam}.
Note from \eqref{4.5:epal} that
\begin{equation}\label{4.23}
  \eps_\B=\alpha_\B+o(\alpha_\B)\ \ \hbox{as}\ \ \B\nearrow\B^*,
\end{equation}
where $\alpha_\B>0$ is defined by \eqref{1.16:alp}.
We then deduce from \eqref{1.16:alp} and \eqref{4.21}--\eqref{4.23} that
\begin{equation}\label{4.24}
\begin{split}
\l(\frac{\alpha_\B}{\eps_\B}\r)^4
=&1+\frac{B_2}{B_1}\eps_\B^2
+\Big(\frac{B_3}{B_1}
+\frac{(1-4\gamma_1\gamma_2)B_1}
      {(2\B^*-a_1-a_2)4\gamma_1^2\gamma_2^2}\Big)\eps_\B^4
+o(\eps_\B^{4})\\
:=&1+\frac{B_2}{B_1}\eps_\B^2+B_4\eps_\B^4+o(\eps_\B^{4})
\ \ \hbox{as}\ \ \B\nearrow\B^*.
\end{split}
\end{equation}

Substituting \eqref{4.23} into \eqref{4.24} yields that
\begin{equation}\label{4.25}
  \l(\frac{\alpha_\B}{\eps_\B}\r)^4=1+\frac{B_{2}}{B_{1}}\alpha_\B^2+o\l(\alpha_\B^2\r)\ \ \hbox{as}\ \ \B\nearrow\B^*,
\end{equation}
which further implies that
\begin{equation}\label{4.26}
  \eps_\B=\alpha_\B-\frac{B_{2}}{4B_{1}}\alpha_\B^3
  +o\l(\alpha_\B^3\r)\ \ \hbox{as}\ \ \B\nearrow\B^*.
\end{equation}
Similar to \eqref{4.25}, one can obtain from \eqref{4.24} and \eqref{4.26} that
\begin{equation*}
\begin{aligned}
  \l(\frac{\alpha_\B}{\eps_\B}\r)^4
  &=1+\frac{B_2}{B_1}\alpha_\B^2+\frac{2B_1^2B_4-B_2^2}{2B_1^2}\alpha_\B^4
+o\l(\alpha_\B^4\r)\ \ \hbox{as}\ \ \B\nearrow\B^*,
\end{aligned}
\end{equation*}
from which we derive that \eqref{4.17} holds true,
and Lemma \ref{lem4.3} is thus proved.
\qed

\vskip 0.1truein

Following Lemmas \ref{lem4.2} and \ref{lem4.3}, we are now going to complete the proof of Theorem \ref{thm3}.

\noindent {\bf Proof of Theorem \ref{thm3}.}
We first obtain from Lemmas \ref{lem4.2} and \ref{lem4.3} that
\begin{equation}\label{4.27}
\begin{aligned}
v_{j\B}(x):=&\sqrt{a^*}\eps_\B u_{j\B}\big(\eps_\B x+x_\B\big)
    e^{-i\big(\frac{\eps_\B\Om}{2}x\cdot x_\B^\perp-\theta_{j\B}\big)}\\
=&\rho_{j\B}Q(x)+\rho_{j\B}\alpha_\B^4\psi_{1}(x)+o(\alpha_\B^4)
\ \ \hbox{as}\ \ \B\nearrow\B^*, \ \ j=1,2,
\end{aligned}
\end{equation}
where $\eps_\B>0$ is defined by \eqref{2.7:eps},  $x_\B=(0,p_\B)\in\R^2$ is the unique maximal point of $|u_{1\B}|^2+|u_{2\B}|^2$, $\theta_{j\B}\in [0,2\pi)$ is chosen such that \eqref{2.10:theta} holds true, $\rho_{j\B}>0$ is as in \eqref{1.23:rhojB}, $\alpha_\B>0$ is defined by \eqref{1.16:alp}, and $\psi_{1}(x)\in C^2(\R^2,\C)\cap L^\infty(\R^2,\C)$ is uniquely given by \eqref{4.11:psik}.

Note from Lemma \ref{lem4.3} that
\begin{equation}\label{4.28}
\frac{\alpha_\B}{\eps_\B}
=1+\frac{B_2}{4B_1}\alpha_\B^2+\frac{8B_1^2B_4-7B_2^2}{32B_1^2}\alpha_\B^4
  +o(\alpha_\B^4)\ \ \hbox{as}\ \ \B\nearrow\B^*.
\end{equation}
We then deduce from \eqref{4.27} and \eqref{4.28} that for $j=1,2$,
\begin{equation}\label{4.29}
\begin{aligned}
&\sqrt{a^*}\alpha_\B u_{j\B}(\alpha_\B x+x_\B)
e^{-i\big(\frac{\alpha_\B\Om}{2}x\cdot x_\B^\perp-\theta_{j\B}\big)}\\
=&\frac{\alpha_\B}{\eps_\B}v_{j\B}\Big(\frac{\alpha_\B}{\eps_\B}x\Big)\\
=&\frac{\alpha_\B}{\eps_\B}\rho_{j\B}\Big(Q\big(\frac{\alpha_\B}{\eps_\B}x\big)
     +\alpha_\B^4\psi_{1}\big(\frac{\alpha_\B}{\eps_\B} x\big)+o(\alpha_\B^4)\Big)\\
=&\Big(1+\frac{B_2}{4B_1}\alpha_\B^2
      +\frac{8B_1^2B_4-7B_2^2}{32B_1^2}\alpha_\B^4+o(\alpha_\B^4)\Big)\\
&\quad\cdot\rho_{j\B}\Big[Q(x)+\Big(\frac{\alpha_\B}{\eps_\B}-1\Big)(x\cdot\nabla Q)
                  +\Big(\frac{\alpha_\B}{\eps_\B}-1\Big)^2\frac{x^T(\nabla^2 Q) x}{2}
                  +\alpha_\B^4\psi_{1}(x)+o(\alpha_\B^4)\Big]\\
=&\rho_{j\B}\Big\{Q(x)
+\alpha_\B^2\frac{B_2}{4B_1}\big(Q+x\cdot\nabla Q\big)\\
&\quad+\alpha_\B^4\Big[\frac{8B_1^2B_4-7B_2^2}{32B_1^2}\big(Q+x\cdot\nabla Q\big)
             +\frac{B_{2}^2}{32B_{1}^2}x^T(\nabla^2 Q)x
          +\frac{B_{2}^2}{16B_{1}^2}\big(x\cdot\nabla Q\big)+\psi_{1}(x)\Big] \Big\}\\
&\quad+o(\alpha_\B^4)\\
:=&\rho_{j\B}\Big(Q(x)+\alpha_\B^2C_{1}(x)+\alpha_\B^4C_{2}(x)\Big)+o(\alpha_\B^4)
\ \ \hbox{as}\ \ \B\nearrow\B^*.
\end{aligned}
\end{equation}

Define for $j=1,2$,
\begin{equation}\label{4.30}
\widetilde{v}_{j\B}(x):=\sqrt{a^*}\alpha_\B u_{j\B}(\alpha_\B x+x_\B)
e^{-i\big(\frac{\alpha_\B\Om}{2}x\cdot x_\B^\perp-\widetilde{\theta}_{j\B}\big)},
\end{equation}
where $\alpha_\B>0$ is defined by \eqref{1.16:alp}, $x_\B=(0,p_\B)$ is the unique maximal point of $|u_{1\B}|^2+|u_{2\B}|^2$, and $\widetilde{\theta}_{j\B}\in[0,2\pi)$ is chosen such that
\begin{equation}\label{4.31:ortho}
\big\|\widetilde{v}_{j\B}-\sqrt{\gamma_j}Q\big\|_{L^2(\R^2)}
=\min\limits_{\theta\in[0,2\pi)}\big\|
\sqrt{a^*}\alpha_\B u_{j\B}(\alpha_\B x+x_\B)
e^{-i\big(\frac{\alpha_\B\Om}{2}x\cdot x_\B^\perp-\theta\big)}
-\sqrt{\gamma_j}Q\big\|_{L^2(\R^2)}.
\end{equation}
By \eqref{4.31:ortho}, we obtain that
\begin{equation}\label{4.32}
    Re\Big(\inte \widetilde{v}_{j\B}(iQ)\Big)=0,\ \ j=1,2.
\end{equation}
Similar to \eqref{2.66}, one can deduce from \eqref{4.27} and \eqref{4.29}--\eqref{4.31:ortho} that
\begin{equation}\label{4.33}
\lim_{\B\nearrow\B^*}|\widetilde{\theta}_{j\B}-\theta_{j\B}|=0,\ \ j=1,2.
\end{equation}
Using \eqref{4.29}, we derive from \eqref{4.30} and \eqref{4.32} that for $j=1,2$,
\begin{equation*}
\begin{aligned}
0&=Re\Big(\inte \widetilde{v}_{j\B}(iQ)\Big)\\
&=Re\Big(\inte\sqrt{a^*}\alpha_\B u_{j\B}(\alpha_\B x+x_\B)
          e^{-i\big(\frac{\alpha_\B\Om}{2}x\cdot x_\B^\perp-\widetilde{\theta}_{j\B}\big)}(iQ)\Big)\\
&=Re\Big(\inte\sqrt{a^*}\alpha_\B u_{j\B}(\alpha_\B x+x_\B)
e^{-i\big(\frac{\alpha_\B\Om}{2}x\cdot x_\B^\perp-\theta_{j\B}\big)}
       e^{-i(\theta_{j\B}-\widetilde{\theta}_{j\B})}(iQ)\Big)\\
&=\rho_{j\B}Re\Big(\inte\Big(Q(x)+\alpha_\B^2C_{1}(x)+\alpha_\B^4C_{2}(x)\Big)
                   e^{-i(\theta_{j\B}-\widetilde{\theta}_{j\B})}(iQ)\Big)
  +o(\alpha_\B^4)\\
&=\rho_{j\B}\sin\big(\theta_{j\B}-\widetilde{\theta}_{j\B}\big)
\inte\Big(Q(x)+\alpha_\B^2C_{1}(x)+\alpha_\B^4C_{2}(x)\Big)Q
   +o(\alpha_\B^4)\ \ \hbox{as}\ \ \B\nearrow\B^*,
\end{aligned}
\end{equation*}
together with \eqref{4.33} then yield that
\begin{equation}\label{4.34}
|\widetilde{\theta}_{j\B}-\theta_{j\B}|=o(\alpha_\B^4)
\ \ \hbox{as}\ \ \B\nearrow\B^*,\ \ j=1,2.
\end{equation}
Combining \eqref{4.29} with \eqref{4.34} then yields that \eqref{1.24:exp} holds true, and the proof of Theorem \ref{thm3} is therefore complete.
\qed

%
%



\begin{thebibliography}{99}
\bibitem{A} A. Aftalion,  Vortices in Bose-Einstein condensates, Progress in Nonlinear Differential Equations and their Applications, 67. Birkh$\ddot{a}$user Boston, Inc., Boston, MA,  2006.



\bibitem{AMW} A. Aftalion, P. Mason and J. C. Wei, {\em Vortex-peak interaction and lattice shape in rotating two-component Bose-Einstein condensates},  Phys. Rev. A {\bf 85} (2012), 033614.



\bibitem{ANS}  A. Aftalion, B. Noris and C. Sourdis, {\em Thomas-Fermi approximation for coexisting two component Bose-Einstein condensates and nonexistence of vortices for small rotation},  Comm. Math. Phys. {\bf 336} (2015), no. 2, 509--579.

\bibitem{AS} A. Aftalion and E. Sandier, {\em Vortex patterns and sheets in segregated two component Bose-Einstein condensates}, Calc. Var. Partial Differ. Eqns. {\bf59} (2020), No. 19, 38 pp.


\bibitem{BC} W. Z. Bao and Y. Y. Cai, {\em Mathematical theory and numerical methods for Bose-Einstein condensation}, Kinet. Relat. Models {\bf6} (2013), 1--135.


\bibitem{BDELL} D. Bonheure, J. Dolbeault,  M. J. Esteban, A. Laptev, M. Loss, {\em Symmetry results in two-dimensional inequalities for Aharonov-Bohm magnetic fields}, Comm. Math. Phys. {\bf 375} (2020),  2071--2087.

\bibitem{BNV} D. Bonheure, M. Nys and J. Van Schaftingen, {\em Properties of ground states of non-linear Schr\"{o}dinger equations under a weak constant magnetic field}, J. Math. Pures Appl. {\bf 124} (2019), 123--168.


\bibitem{BSSD} V. Bretin, S. Stock, Y. Seurin and J. Dalibard, {\em Fast rotation of a Bose-Einstein condensate}, Phys. Rev. Lett. {\bf 92} (2004), 050403.


\bibitem{C} T. Cazenave, Semilinear Schr\"{o}dinger Equations, Courant Lecture Notes in Mathematics  Vol. 10, Courant Institute of Mathematical Science/AMS, New York, 2003.



\bibitem{Cing} S. Cingolani and S. Secchi, {\em Semiclassical limit for nonlinear Schr\"{o}dinger equations with electromagnetic fields}, J. Math. Anal. Appl.
{\bf 275} (2002), 108--130.

\bibitem{Cing2}  S. Cingolani and S. Secchi, {\em Semiclassical states for NLS equations with magnetic potentials having polynomial growths}, J. Math. Phys.
{\bf 46} (2005), 053503.


\bibitem{DW} E. N. Dancer and J. C. Wei, {\em Spike solutions in coupled nonlinear Schr\"{o}dinger equations with attractive interaction}, Trans. Amer. Math. Soc. {\bf361} (2009),  1189--1208.



\bibitem{F} A. L. Fetter, {\em Rotating vortex lattice in a Bose-Einstein condensate trapped in combined
quadratic and quartic radial potentials}, Phys. Rev. A {\bf 64} (2001), 063068.


\bibitem{F2} A. L. Fetter, B. Jackson and S. Stringari, {\em Rapid rotation of a Bose-Einstein condensate in
a harmonic plus quartic trap}, Phys. Rev. A {\bf 71} (2005), 013605.

\bibitem{GGLL}Y. S. Gao, Y. J. Guo, Y. Li and Y. Luo, {\em Nonexistence of Vortices for Rotating Two-Component Focusing Bose Gases}, submitted. (2022), 59 pages, arXiv:2211.14808.


\bibitem{GRPG}  J. J. Garc{\'{\i}}a-Ripoll, V. M. P\'{e}rez-Garc{\'{\i}}a and F. Sols, {\em Split vortices in optically coupled Bose-Einstein condensates},  Phys. Rev. A {\bf66} (2002), 021602.


\bibitem{GNN} B. Gidas, W. M. Ni and L. Nirenberg, {\em Symmetry of positive solutions of nonlinear elliptic equations in $\R^{n}$}, Mathematical analysis and applications Part A, Adv. Math. Suppl. Stud. Vol. 7, Academic Press, New York (1981), 369--402.


\bibitem{G} M. Grossi, {\em On the number of single-peak solutions of the nonlinear Schr\"{o}dinger equations}, Ann. Inst. H. Poincar\'{e} Anal. Non Lin\'{e}aire {\bf19} (2002), 261--280.

\bibitem{GX} Q. Guo and H. F. Xie, {\em Existence and local uniqueness of normalized solutions for two-component Bose-Einstein condensates},  Z. Angew. Math. Phys. {\bf72} (2021), 189.




\bibitem{GLLP} Y. J. Guo, Y. Li, Y. Luo and S. J. Peng, {\em Axial Symmetry of Normalized Solutions for Magnetic Gross-Pitaevskii Equations with Anharmonic Potentials}, submitted (2023), 40 pages, arXiv:2310.00556.

\bibitem{GLWZ1} Y. J. Guo, S. Li, J. C. Wei and X. Y. Zeng, {\em Ground states of two-component attractive Bose-Einstein condenstates I: Existence and uniqueness}, J. Funct. Anal. {\bf276} (2019), 183--230.


\bibitem{GLW} Y. J. Guo, C. S. Lin and J. C. Wei, {\em Local uniqueness and refined spike profiles of ground states for two-dimensional attractive Bose-Einstein condensates}, SIAM J. Math. Anal. {\bf49} (2017), 3671--3715.

\bibitem{GLP} Y. J. Guo, Y. Luo and S. J. Peng, {\em Local uniqueness of ground states for rotating Bose-Einstein condenstates with attractive interactions}, Calc. Var. Partial Differ. Eqns. {\bf60} (2021), 237.

\bibitem{GLY} Y. J. Guo, Y. Luo and W. Yang, {\em The nonexistence of vortices for rotating Bose-Einstein condenstates with attractive interactions},  Arch. Ration. Mech. Anal. {\bf238} (2020), 1231--1281.



\bibitem{GZZ1} Y. J. Guo, X. Y. Zeng and H. S. Zhou, {\em Energy estimates and symmetry breaking in attractive Bose-Einstein condensates with ring-shaped potentials}, Ann. Inst. H. Poincar\'{e} Anal. Non Lin\'{e}aire {\bf33}, (2016), 809--828.









\bibitem{KMKR} E. \"{O}. Karabulut, F. Malet, G. M. Kavoulakis and S. M. Reimann, {\em Phase diagram of a rapidly rotating two-component Bose gas}, Phys. Rev. A {\bf87} (2013), 043609.


\bibitem{KTU1} K. Kasamatsu, M. Tsubota and M. Ueda, {\em Vortex phase diagram in rotating two-component Bose-Einstein condensates}, Phys. Rev. Lett. {\bf91} (2003), 150406.




\bibitem{K} M. K. Kwong, {\em Uniqueness of positive solutions of $\Delta u-u+u^{p}=0$ in $\R^{N}$},  Arch. Ration. Mech. Anal. {\bf105} (1989), 243--266.



\bibitem{LS} E. H. Lieb and R. Seiringer, {\em Derivation of the Gross-Pitaevskii equation for rotating Bose gases}, Comm. Math. Phys. {\bf264} (2006), 505--537.

\bibitem{LS1} E. H. Lieb and J. P. Solovej, {\em Ground state energy of the two-component charged Bose gas}, Comm. Math. Phys. {\bf252} (2004), 485--534.


\bibitem{LW}  T. C. Lin and J. C. Wei, {\em Spikes in two-component systems of nonlinear Schr\"{o}dinger equations with trapping potentials},  J. Differ. Eqns. {\bf229} (2006), 538--569.


\bibitem{LZ} Y. Liu and S. Y. Zhang, {\em The ground states and pseudospin textures of rotating two-component Bose-Einstein condensates trapped in harmonic plus quartic potential}, Chin. Phys. B {\bf 25} (2016), 090304.


\bibitem{LPWY} P. Luo, S. J. Peng, J. C. Wei and S. S. Yan, {\em Excited states of Bose-Einstein condensates with degenerate attractive interactions}, Calc. Var. Partial Differ. Eqns. {\bf 60} (2021), Paper No. 155, 33 pp.









\bibitem{R} J. Royo-Letelier, {\em Segregation and symmetry breaking of strongly coupled two component Bose-Einstein condensates in a harmonic trap}, Calc. Var. Partial Differ. Eqns. {\bf49} (2014), 103--124.




%
%


\bibitem{SBCD}  S. Stock, V. Bretin,  F. Chevy and J. Dalibard, {\em Shape oscillation of a rotating Bose-Einstein condensate}, Europhys. Lett. {\bf 65}  (2004), 594--600.

\bibitem{TKT} H. Takeuchi, K. Kasamatsu and M. Tsubota, {\em Vortex Structures in Rotating Two-Component Bose-Einstein Condensates in an Anharmonic Trapping Potential}, AIP Conference Proceedings {\bf850}, (2006),  57--58.


\bibitem{WY} J. C. Wei and W. Yao, {\em Uniqueness of positive solutions to some coupled nonlinear Schr\"{o}dinger equations}, Commun. Pure Appl. Anal. {\bf11} (2012), 1003--1011.

\bibitem{ZZ} X. Zeng and H. Zhou, {\em Uniqueness of single peak solutions for coupled nonlinear Gross-Pitaevskii equations with potentials}, submitted (2022), 31 pages, arXiv:2204.11711.



\end{thebibliography}
\end{document}